\documentclass[a4paper,reqno]{amsart}%
\usepackage{latexsym}
\usepackage{amsmath}
\usepackage{graphicx}
\usepackage{amsfonts}%
\usepackage{amssymb}
\theoremstyle{plain}
\begingroup
\newtheorem{theorem}{Theorem}[section]

\newtheorem{proposition}[theorem]{Proposition}
\newtheorem{corollary}[theorem]{Corollary}
\endgroup
\theoremstyle{definition}
\begingroup
\newtheorem{definition}[theorem]{Definition}
\newtheorem{remark}[theorem]{Remark}

\endgroup
\theoremstyle{remark}
\begingroup
\endgroup
\mathchardef\emptyset="001F
\numberwithin{equation}{section}
\newcommand{\e}{\varepsilon}
\newcommand{\Om}{\Omega}
\newcommand{\R}{{\mathbb {R}}}

\newcommand{\N}{\mathbb {N}}

\newcommand{\supp}{{\rm supp}}
\newcommand{\Lone}{{\mathcal {L}}^1}

\newcommand{\LM}[2]{\text{\vrule width.4pt \vbox to#1pt{\vfill
\hrule width#2pt height.4pt}}}
\newcommand{\LLL}{{\mathchoice {\>\LM{7}{5}\>}{\>\LM{7}{5}\>}{\,\LM{5}{3.5}\,}{\,\LM{3.35}{2.5}\,}}}
\title[A higher order model for image restoration: the one dimensional case]
{A higher order model for image restoration:\\ the one dimensional case}
\author{G.\ Dal Maso}
 \address[Gianni Dal Maso]{SISSA, Via Beirut 2,
34014 Trieste, Italy}
 \email[Gianni Dal Maso]{dalmaso@sissa.it}
\author{I.\ Fonseca}
\address[Irene Fonseca]{Carnegie Mellon University, Pittsburgh PA 15213-3890, USA}
\email[Irene Fonseca]{fonseca@andrew.cmu.edu}
\author{G.\ Leoni}
\address[Giovanni Leoni]{Carnegie Mellon University, Pittsburgh PA 15213-3890, USA}
\email[Giovanni Leoni]{giovanni@andrew.cmu.edu}
\author{M.\ Morini}
 \address[Massimiliano Morini]{SISSA, Via Beirut 2,
34014 Trieste, Italy}
\email[Massimiliano Morini]{morini@sissa.it}
\begin{document}

\begin{abstract}
The one-dimensional version of the higher order total
variation-based model for image restoration proposed by Chan,
Marquina, and Mulet in \cite{Cha-Mar-Mul} is analyzed. A suitable
functional framework in which the minimization problem is well posed
is being proposed and it is proved analytically that the higher
order regularizing term prevents the occurrence of the {\em
staircase effect}.
\end{abstract}
\maketitle

{\small

\bigskip

}

\tableofcontents
\bigskip\bigskip

\section{Introduction}
Deblurring and denoising of images are fundamental problems in image processing and gave rise in the past few years to a vast variety  of techniques and methods touching different fields of mathematics. Among them, variational methods based on the minimization of some energy functional
have been successfully employed to treat a fairly general class of image restoration problems. Typically, such functionals
 present a {\em fidelity term}, which penalizes the distance between the reconstructed image $u$ and the noisy image
 $g$
with respect to a suitable metric, and a regularizing term, which makes high frequency noise energetically
unfavorable.

When the fidelity term is given  by the squared $L^2$ distance multiplied by a parameter $\lambda>0$ and the
regularizing term is represented by the total variation, we are led to the following minimization problem
\begin{equation}\label{ROF0}
\min\Big\{|Du|(\Om)+\lambda\int_{\Om}|u-g|^2\,dx:\, u\in BV(\Om)\Big\}\,,
\end{equation}
which was proposed by Rudin, Osher, and Fatemi in \cite{ROF}. Here $\Om$ is an open bounded  domain in one or two dimensions,
$BV(\Om)$ denotes the space of functions of bounded variations in $\Om$, and $|Du|(\Om)$ stands for the total
variation of $u$ in $\Om$.  The main feature of the total variation-based image restoration is perhaps represented
by the tendency to yield (almost) piecewise-constant solutions or, in other words, ``blocky'' images. Typically one observes that {\em ramps} (i.e., affine regions) in the original image give rise to  staircase-like structures in the reconstructed
 image, a  phenomenon which is often referred to
 as the {\em staircase effect}. This means that the original edges are well preserved by this method, but also that
 many artificial discontinuities can be generated by the presence of noise, while the finer details of the
 objects contained in the image may not be properly recovered.

Several variants of \eqref{ROF0} have been subsequently proposed in order to fix these drawbacks. In this paper
we follow the approach of Chan, Marquina, and Mulet \cite{Cha-Mar-Mul}: Since the total variation does not
distinguish between jumps and smooth transitions their idea is to consider an additional  penalization of the discontinuities  by taking second derivatives into account. More precisely, they propose a regularizing term of the form
\begin{equation}\label{regtermCMM}
 \int_{\Om}|\nabla u|\, dx+ \int_{\Om}\psi(|\nabla u|)h(\Delta u)\, dx\,,
\end{equation}
where $\psi$ is a function that must satisfy suitable conditions at infinity
in order to allow  jumps.

In this paper we consider the following 1-D version of \eqref{regtermCMM}:
\begin{equation}\label{1dCMM}
\mathcal {F}_p(u):=\int_a^b|u'|\, dx+\int_a^b\psi(|u'|)|u''|^p\, dx\,,
\end{equation}
where $a<b$ are real numbers  and  $p\in [1,+\infty)$.
Our main analytical objective is twofold:
\begin{itemize}
\item[(i)] to set up a proper functional framework where the minimization
problem corresponding to
\begin{equation*}
\mathcal {F}_p(u)+\lambda\int_{\Om}|u-g|^2\,dx
\end{equation*}
is well posed;
\item[(ii)] to give an analytical proof of the fact the  higher order regularizing term
eliminates the staircase effect.
\end{itemize}
We point out  here that we carry out the first part of this program
by using the theory of relaxation (see \cite{DM} for a general
introduction): We regard $\mathcal {F}_p$ as defined
 for all functions in the Sobolev space $W^{2,p}({]a,b[})$, we extend it to $L^1({]a,b[})$ by setting
$\mathcal {F}_p(u):=+\infty$ if $u\in  L^1({]a,b[})\setminus W^{2,p}({]a,b[})$, and then we identify its lower
 semicontinuous envelope with respect to the strong $L^1$ convergence. The extension of our results
to higher dimensions will be the subject of a subsequent paper.

For completeness we conclude by mentioning that  other approaches
have been considered to avoid staircasing: The works
 by Geman and Reynolds \cite{GR} and Chambolle and Lions \cite{CL} contain a different use of higher order
derivatives as regularizing terms; in \cite{BCM}, Blomgren, Chan, and Mulet propose a $BV$-$H^1$ interpolation approach, while Kindermann, Osher, and Jones avoid in \cite{KOJ} the use of second derivatives by considering a
sort of nonlocal total variation.

The plan of the paper is the following: In Section~\ref{s:p=1} we
consider the case $p=1$; i.e., we identify the relaxation
$\overline{\mathcal{F}}_1$ of $\mathcal{F}_1$, while in
Section~\ref{s:p>1} we deal with the case $p>1$. The analysis turns
out to be considerably more delicate in the former case. Moreover,
the domains of the relaxed
 functionals are quite peculiar (see Definitions~\ref{def space} and \ref{def spacep>1}) and display
properties which are qualitatively different in the two cases. In
particular, it turns out that piecewise constant functions
corresponding to images with genuine edges are approximable by
sequences with bounded energy only for $p=1$. Finally, in
Section~\ref{s:staircase} we investigate the staircase effect. After
exhibiting an analytical example of staircasing for the
Rudin-Osher-Fatemi model (Theorem~\ref{th:ROFstaircase}), we prove
that the new model does indeed prevent the occurrence of this
phenomenon. More precisely, we show that
 whenever the datum $g$ is of the form $g=g_1+h$, with $g_1$ a regular image and $h$ a highly oscillating
 noise, the
 reconstructed  image is regular as well (Theorems~\ref{th:m1} and \ref{theorem non staircase p>1}).

\section{The case $p=1$}
\label{s:p=1}

We start by  studying the compactness properties and the relaxation
of \eqref{1dCMM} in the case $p=1$. Throughout this section
$\psi\colon{\mathbb{R}}\rightarrow{]0,+\infty\lbrack}$ will be a
bounded Borel function such that
\begin{equation}
M:=\int_{-\infty}^{+\infty}\psi(t)\,dt<+\infty\label{intfin}%
\end{equation}
and
\begin{equation}
\inf_{t\in K}\psi(t)>0\quad\text{for every compact set }K\subset{\mathbb{R}%
}\,. \label{inf>0}%
\end{equation}
Let $\Psi_1\colon\overline{{\mathbb{R}}}\rightarrow [0,M]$ be the
increasing function defined by
\begin{equation*}
\Psi_1(t):=\int_{-\infty}^{t}\psi(s)\,ds%
\end{equation*}
and let $\Psi^{-1}_1\colon [0,M]\rightarrow\overline{{\mathbb{R}}}$ be its
inverse function.

Given  a bounded open interval ${]a,b[}$ in ${\mathbb{R}}$, we consider the
functional ${\mathcal{F}}_1\colon L^{1}({]a,b[})\rightarrow[0,+\infty]$
defined by
\begin{equation}
{\mathcal{F}}_1(u):=%
\begin{cases}
\displaystyle  \int_{a}^{b}|u^{\prime}|\,dx+\int_{a}^{b}\psi(u^{\prime
})|u^{\prime\prime}|\,dx & \text{if }u\in W^{2,1}({]a,b[})\,,\\
+\infty & \text{otherwise.}%
\end{cases}
\label{functional}%
\end{equation}

The first step in the study of \eqref{functional} will consist in identifying  the subspace of $L^1$ functions   which can be approximated by energy bounded sequences. In order to do so we need to introduce some notation and recall some basic facts about $BV$ functions of one variable. This will be the content of the next subsection.

\subsection{$BV$ functions of one variable}
  We recall that a function $u\in L^1({]a,b[})$ belongs to $BV({]a,b[})$ if and only if
\begin{equation}\label{essTV}
\sup\Big\{\int_a^b u\varphi'\, dx:\, \varphi\in C^1_c(]a,b[)\,,\, |\varphi|\leq 1\Big\}<+\infty\,.
\end{equation}
Note that this implies that the distributional derivative $u'$ of
$u$ is a bounded Radon measure in ${]a,b[}$. We will often consider
the Lebesgue decomposition
$$
u^{\prime}=(u^{\prime})^{a}\Lone+(u^{\prime})^{s}
$$
where $(u^{\prime})^{a}$ is the density of the absolutely continuous part of $u'$ with respect to the Lebesgue
 measure $\Lone$ on ${]a,b[}$, while $(u^{\prime})^{s}$ is its singular part. We will denote  the total variation
measure of $u'$ by $|u'|$. In particular,  $|u'|({]a,b[})$ equals the value of the supremum in \eqref{essTV}.
For every function $u\in BV({]a,b[})$ the following left and right approximate limits
$$
u_-(y):=\lim_{\varepsilon\rightarrow0+}\frac{1}{\varepsilon
}\int_{y-\varepsilon}^{y}u(x)\,dx\,,\qquad
 u_+(y):=\lim_{\varepsilon\rightarrow 0+}\frac{1}{\varepsilon}\int
_{y}^{y+\varepsilon}u(x)\,dx
$$
are well defined at every point $y\in{]a,b[}$. In fact, $u_-(y)$ is well defined also at $y=b$ while
$u_+(y)$ exists also at $y=a$. The functions $u_-$ and $u_+$ coincide $\Lone$-a.e.\ with $u$ and are
 left and right continuous, respectively. Moreover, it turns out that the set $S_u:=\{y\in{]a,b[}:\, u_-(y)\neq
 u_+(y)\}$ is at most countable. The set $S_u$ is often referred to as the set of essential discontinuities or
{\em jump points} of~$u$.

It is well known that, in turn, the singular part $(u')^s$ splits into the sum of an atomic measure concentrated
on   $S_u$ and a singular diffuse measure $(u^{\prime})^{c}$, called the {\em Cantor part} of $u'$:
$$
(u^{\prime})^{s}=[u]{\mathcal{H}}^{0}
\LLL S_{u}+(u^{\prime})^{c}\,,
$$
where we set $[u]:=u_+-u_-$ and ${\mathcal{H}}^{0}$ stands for the counting measure. Finally, we recall that every $u\in BV({]a,b[})$ is differentiable at $\Lone$-a.e.\ $y$ in
${]a,b[}$ with derivative given by $(u^{\prime})^{a}(y)$. In this case we will often write, with a slight abuse of notation, $u^{\prime}(y)$ instead of $(u^{\prime})^{a}(y)$.
\par
We say that a sequence $\{u_{k}\}$ of functions in $BV({]a,b[})$
{\em weakly star converges in } $BV({]a,b[})$ to a function $u\in
BV({]a,b[})$ if $u_n\to u$ in $L^1({]a,b[})$ and $u'_k
\rightharpoonup u'$ weakly$^*$ in $M_b({]a,b[})$, where
$M_b({]a,b[})$ is the space of bounded Radon measures.
\par
We will also need sometimes the notion of total variation for a function defined everywhere. We recall that
$u\colon{]a,b[}\to\R$ has bounded {\em  pointwise total variation} over the interval
${]c,d[}\subset{]a,b[}$ if
$$
{\rm Var\,} (u; ]c,d[):= \sup \sum_{i=1}^k |u(y_i)-u(y_{i-1})|<+\infty\,,
$$
where the supremum is taken over all finite families
$y_0,y_1,\dots,y_k$ such that $c<y_0<y_1<\dots<y_k<d$, $k\in\N$. It
is easy to see that if $u$ has bounded  pointwise total variation in
${]a,b[}$, then it admits left and right limits at every point, it
belongs to $BV({]a,b[})$, and $|u'|(]c,d[)\leq {\rm Var\,} (u;
]c,d[)$ for every interval ${]c,d[}\subset{]a,b[}$. Conversely, if
$u\in BV({]a,b[})$, the {\em precise representatives} $u_-$ and
$u_+$ have bounded pointwise total variation and satisfy
$$
|u'|(]c,d[)={\rm Var\,} (u_-; ]c,d[)={\rm Var\,} (u_+; ]c,d[)
$$
for every interval ${]c,d[}\subset{]a,b[}$.\par Finally, we recall
the Helly theorem: For every  bounded sequence of functions
$u_k:{]a,b[}\to\R$ such that $\sup_k{\rm Var\,} (u_k;
{]a,b[})<+\infty$, there exist $u$, with pointwise total variation
in ${]a,b[}$, and a subsequence (not relabeled) such that $u_k\to u$
pointwise.

We refer to \cite{Ru} and \cite{HS} for an exhaustive exposition of the
properties of $BV$ functions of one variable.

\subsection{Compactness}
To define the subspace of $L^1$ functions that can be approximated by energy bounded sequences, for every function $u\in BV\left(  \left]  a,b\right[
\right)$ we consider  the sets
\begin{align}
&  \displaystyle  Z^{+}[(u^{\prime})^{a}]:=\Big\{  x\in{]a,b[}:
\lim_{\varepsilon\rightarrow0+}\frac{1}{2\varepsilon}
\int_{x-\varepsilon}^{x+\varepsilon}(u^{\prime})^{a}\,dx=+\infty\Big\}  \,,\label{Zu}\\
&  \displaystyle  Z^{-}[(u^{\prime})^{a}]:=\Big\{  x\in{]a,b[}:
\lim_{\varepsilon\rightarrow0+}\frac{1}{2\varepsilon}\int_{x-\varepsilon
}^{x+\varepsilon}(u^{\prime})^{a}\,dx=-\infty\Big\}  \,. \label{Zu-}%
\end{align}
It is also convenient to define%
\begin{equation*}
Z[(u^{\prime})^{a}]:=Z^{+}[(u^{\prime})^{a}]\cup
Z^{-}[(u^{\prime})^{a}] \,.
\end{equation*}

\begin{definition}
\label{def space}Let $X^1_{\psi}({]a,b[})$ be the set of all functions $u\in
BV({]a,b[})$ such that $v:=\Psi_1\circ (u^{\prime})^{a}$ belongs to $BV({]a,b[})$
and the positive part $\left(  (u^{\prime})^{c}\right)^{+}$ and the negative
part $\left(  (u^{\prime})^{c}\right)^{-}$ of the measure $(u^{\prime})^{c}$
are concentrated on $Z^{+}[(u^{\prime})^{a}]$ and $Z^{-}[(u^{\prime})^{a}]$, respectively.
\end{definition}

\begin{remark}
\label{remark limits u'}Note that if $u\in X^1_{\psi}({]a,b[})\,$then the
limits
\begin{equation}
(u^{\prime})_{-}^{a}(y):=\lim_{\varepsilon\rightarrow0+}\frac{1}{\varepsilon
}\int_{y-\varepsilon}^{y}(u^{\prime})^{a}\,dx\,,\qquad(u^{\prime})_{+}%
^{a}(y):=\lim_{\varepsilon\rightarrow0+}\frac{1}{\varepsilon}\int
_{y}^{y+\varepsilon}(u^{\prime})^{a}\,dx\,\label{uprime+}%
\end{equation}
exist in $\overline{\mathbb{R}}$ for every $y$. More precisely, $(u^{\prime})_{-}^{a}$ exists also at
 $y=b$ while $(u^{\prime})_{+}^{a}$ is well defined also at $y=a$.
 Indeed, since $v=\Psi_1\circ
(u^{\prime})^{a}$ is a $BV$ function, it admits a precise
representative $\tilde{v}$ such that the right and left limits exist
at every point, and the same property holds for
$\Psi^{-1}_1(\tilde{v})$. As $\Psi^{-1}_1(\tilde
{v})=(u^{\prime})^{a}$ ${\mathcal{L}}^{1}$-a.e.\ in ${]a,b[}$, the
limits considered in \eqref{uprime+} are everywhere well-defined.
Moreover the set $S_{(u^{\prime})^{a}}:=S_{v}$ is at most countable
and
\begin{equation}\label{12300}
(u^{\prime})_{-}^{a}=(u^{\prime})_{+}^{a}\text{\quad on }
{]a,b[}\setminus
S_{(u^{\prime})^{a}}.
\end{equation}
We also remark that  $(u^{\prime})_{-}^{a}$ and $(u^{\prime})_{+}^{a}$ are left and right continuous, which, in turn,  implies that the functions defined by
$$
\begin{array}{c}
 (u^{\prime})_{\vee}^{a}(x){:=}\max\left\{  (u^{\prime})_{+}^{a}(x),(u^{\prime})_{-}^{a}(x)\right\}\,,\
  (u^{\prime})_{\wedge}^{a}(x){:=}\min\left\{
(u^{\prime})_{+}^{a}(x),(u^{\prime})_{-}^{a}(x)\right\}\text{ if $x\in{]a,b[}$,}\smallskip\\
 (u^{\prime})_{\vee}^{a}(a)=(u^{\prime})_{\wedge}^{a}(a):=(u^{\prime})_{+}^{a}(a)\,,\quad\text{and }
(u^{\prime})_{\vee}^{a}(b)=(u^{\prime})_{\wedge}^{a}(b):=(u^{\prime})_{-}^{a}(b)
\end{array}
$$
are upper and lower semicontinuous in $[a,b]$, respectively. By
\eqref{12300} we have
\begin{eqnarray*}
Z^{+}[(u^{\prime})^{a}]\setminus S_{(u^{\prime})^{a}}
=\{x\in{]a,b[}:   (u^{\prime})_{\wedge}^{a}(x)=+\infty\} \setminus S_{(u^{\prime})^{a}}\,,
\\
Z^{-}[(u^{\prime})^{a}]\setminus S_{(u^{\prime})^{a}}
=\{x\in{]a,b[}:   (u^{\prime})_{\vee}^{a}(x)=-\infty\} \setminus S_{(u^{\prime})^{a}}\,.
\end{eqnarray*}
Therefore $\left(  (u^{\prime})^{c}\right)^{+}$ is concentrated on the set
$\{x\in{]a,b[}:   (u^{\prime})_{\wedge}^{a}(x)=+\infty\}$ and
$\left(  (u^{\prime})^{c}\right)^{-}$ is concentrated on the set
$\{x\in{]a,b[}:   (u^{\prime})_{\vee}^{a}(x)=-\infty\}$.
\end{remark}

Before we proceed we show that
the space $X^1_{\psi}({]a,b[})$ contains functions with nontrivial Cantor part
when $\psi$ satisfies suitable decay estimates  at infinity.

\begin{proposition}
\label{prop:cantor} Assume that $\psi\colon{\mathbb{R}}\rightarrow
{]0,+\infty\lbrack}$ is a bounded Borel function satisfying $\left(
\ref{intfin}\right)  $, $\left(  \ref{inf>0}\right)  $, and
\begin{equation}
\psi(t)\leq\frac{c}{t^{\alpha}}\label{algebraic}%
\end{equation}
for all $t\geq1$ and for some $c>0$, $\alpha>1$. Then there exists
$u\in X^1_{\psi}({]a,b[})$ with $(u^{\prime})^{c}\neq0$.
\end{proposition}

\begin{proof}
For simplicity we take ${]a,b[}={]0,1[}$.

\textbf{Step 1}: We start by recalling the definition of the generalized
Cantor set $\mathbb{D}_{\delta}$, where $\delta\in{]0,\frac{1}{2}[}$(see for
instance \cite[Chapter 1, Section 2.4]{DiBe}). The construction is
entirely similar to the one of the (ternary) Cantor set with the only
difference that the middle intervals removed at each step have length
$1-2\delta$ times the length of the intervals remaining from the previous
step. To be more precise, remove from $\left[  0,1\right]  $ the interval
$I_{11}:=\left(  \delta,1-\delta\right)  $. At the second step remove from
each of the remaining closed intervals $[0,\delta]$ and $[1-\delta,1]$ the
middle intervals, denoted by $I_{12}$ and $I_{22}$, of length $\delta\left(
1-2\delta\right)  $. Continuing in this fashion at each step $n$ we remove
$2^{n-1}$ middle intervals $I_{1n}$, \ldots, $I_{2^{n-1}n},$ each of length
$\delta^{n-1}\left(  1-2\delta\right)  $. The \emph{generalized} \emph{Cantor
set} $\mathbb{D}_{\delta}$ is defined as%
\[
\mathbb{D}_{\delta}:=\left[  0,1\right]
\setminus\bigcup_{n=1}^{\infty }\bigcup_{k=1}^{2^{n-1}}I_{kn}\,.
\]
The set $\mathbb{D}_{\delta}$ is closed (since its complement is given by a
family of open intervals) and
\[
\mathcal{L}^{1}\left(  \mathbb{D}_{\delta}\right)  =1-\sum_{n=1}^{\infty}%
\sum_{k=1}^{2^{n-1}}\Lone\left(  I_{kn}\right)  =1-\sum
_{n=1}^{\infty}\sum_{k=1}^{2^{n-1}}\delta^{n-1}\left(
1-2\delta\right) =1-\left(  1-2\delta\right)
\sum_{n=1}^{\infty}\left(  2\delta\right) ^{n-1}=0\,.
\]
 Next we recall the
definition of the corresponding Cantor function $f_{\delta}$. Set
\[
g_{n}:=\frac{1}{(2\delta)^{n}}\Big(  1-\sum_{j=1}^{n}\sum_{k=1}^{2^{j-1}}%
\chi_{I_{kj}}\Big) \, ,
\]
and define $f_{n}(x):=\int_{0}^{x}g_{n}(t)\,dt$. It can be shown
that $\{f_{n}\}$ converges uniformly to a continuous nondecreasing
function $f_{\delta}$ such that $f_{\delta}(0)=0$,
$f_{\delta}(1)=1$, and $f_{\delta}^{\prime}=\left(
f_{\delta}^{\prime}\right)  ^{c}$ is supported on
$\mathbb{D}_{\delta}$.

\textbf{Step 2}: We claim that it is enough to find a constant $\delta
\in{]0,\frac{1}{2}[}$ for which it is possible to construct a continuous integrable function $w_{\delta
}:{]0,1[}\rightarrow\lbrack0,+\infty]$ such that $\Psi_1\circ w_{\delta}\in
BV({]0,1[})$ and $w_{\delta}(x)=+\infty$ if and only if $x\in\mathbb{D}%
_{\delta}$. Indeed, setting $u_{\delta}(x):=\int_{0}^{x}w_{\delta
}(t)\,dt+f_{\delta}(x)$, we have that $u_{\delta}\in BV({]0,1[})$, $u_{\delta
}$ is continuous, $\left(  u_{\delta}^{\prime}\right)  ^{a}=w_{\delta}$ so that
$Z^{+}[(u_{\delta}^{\prime})^{a}]=Z[(u_{\delta}^{\prime})^{a}]=\mathbb{D}%
_{\delta}$ and $\Psi_1\circ\left(  u_{\delta}^{\prime}\right)  ^{a}\in
BV({]0,1[})$. Moreover, $(u_{\delta}^{\prime})^{c}=\left(  f_{\delta}^{\prime
}\right)  ^{c}$ is supported on $\mathbb{D}_{\delta}=Z^{+}[(u_{\delta}%
^{\prime})^{a}]$. Hence $u_{\delta}$ belongs to $X^1_{\psi}({]a,b[})$.

\textbf{Step 3}: It remains to construct $w_{\delta}$ for a suitable $\delta
\in{]0,\frac{1}{2}[}$. Consider a convex
function $\phi:{]0,1[}\rightarrow\lbrack0,+\infty)$ such that
\begin{equation}
\lim_{x\rightarrow0^{+}}\phi(x)=\lim_{x\rightarrow1^{-}}\phi(x)=+\infty
,\text{\quad}\phi(\tfrac{1}{2})=0\,, \label{limiti}%
\end{equation}
and
\begin{equation}
\int_{0}^{1}\phi(x)\,dx=1\,. \label{intuno}%
\end{equation}
Choose $s>0$ so large  that
\begin{equation}\label{chooset}
\alpha>\frac {s+1}{s}\,.
\end{equation}
For $x\in I_{kn}$ (see Step 1) define
\begin{equation}
\phi_{kn}(x):=2^{sn}+\phi\left(  \tfrac{x-a_{kn}}{\delta^{n-1}\left(
1-2\delta\right)  }+\tfrac{1}{2}\right) \, , \label{phikn}%
\end{equation}
where $a_{kn}$ is the mid point of the interval $I_{kn}$. Finally set
\[
w_{\delta}:=\sum_{n=1}^{\infty}\sum_{k=1}^{2^{n-1}}\phi_{kn}\chi_{I_{kn}%
}+I_{\mathbb{D}_{\delta}}\,,
\]
where $I_{\mathbb{D}_{\delta}}$ is the indicator function of the set
$\mathbb{D}_{\delta}$, that is,
\[
I_{\mathbb{D}_{\delta}}(x) :=\left\{
\begin{array}
[c]{ll}%
+\infty& \text{if }x\in\mathbb{D}_{\delta}\,,\\
0 & \text{otherwise.}%
\end{array}
\right.
\]
Using
the fact that%
\[
\int_{I_{kn}}\phi_{kn}\,dx=\left(  2^{sn}+1\right)
\delta^{n-1}\left( 1-2\delta\right)\,  ,
\]
which follows from (\ref{intuno}) and a change of variables, we have
\[
\int_{0}^{1}w_{\delta}\,dx=\sum_{n=1}^{\infty}\sum_{k=1}^{2^{n-1}}\left(
2^{sn}+1\right)  \delta^{n-1}\left(  1-2\delta\right)  <\infty
\]
for $\delta<\frac{1}{2^{s+1}}$. To estimate the total variation of $v:=\Psi_1\circ w_{\delta}$ we consider the approximating sequence
\[
v_{m}(x):=\left\{
\begin{array}
[c]{ll}%
\Psi_1\circ\phi_{kn}(x) & \text{if }x\in I_{kn}\text{, }1\leq k\leq2^{n-1}\text{, }
1\leq n\leq m\text{,}\\
M & \text{otherwise.}%
\end{array}
\right.
\]
By \eqref{algebraic}, (\ref{limiti}), (\ref{phikn}), and the convexity of $\phi$ it can be seen that
$$
{\rm Var\,}(v_m; I_{kn} )=2(M-\Psi_1(2^{sn}))=2\int_{2^{sn}}^{+\infty}\psi(t)\, dt\leq \frac{2c}{\alpha-1}\frac{1}{2^{sn(\alpha-1)}}\,.
$$
It follows that
$$
{\rm Var\,}(v_m; ]0,1[)\leq  \frac{2c}{\alpha-1}\sum_{n=1}^{m}\sum_{k=1}^{2^{n-1}}
\frac{1}{2^{sn(\alpha-1)}}\leq   \frac{2c}{\alpha-1}\sum_{n=1}^{\infty}\frac{1}{2^{sn(\alpha-1)-n+1}}\,.
$$
The last series is finite thanks to \eqref{chooset}. Therefore the
$v_m$'s have equibounded total variations and, since $v_m\to v$ in
$L^1(]0,1[)$, we conclude that $v\in BV(]0,1[)$.
\end{proof}

Energy bounded sequences are compact in $X^1_{\psi}({]a,b[})$, as made precise by the following theorem.

\begin{theorem}
\label{theorem compactness} Let $\{u_{k}\}$ be a sequence of
functions bounded in $L^{1}({]a,b[})$ such that
\begin{equation}
C:=\sup_{k}{\mathcal{F}}_1(u_{k})<+\infty\,.\label{supenergy}%
\end{equation}
Then there exist a subsequence (not relabeled) $\{u_{k}\}$ and a
function $u\in X^1_{\psi
}({]a,b[})$ such that%
\begin{align}
& u_{k}\rightharpoonup u\quad\text{weakly}^{\ast}\text{ in }BV({]a,b[}%
)\,,\label{convw*}\\
& \Psi_1\circ u'_k\rightharpoonup \Psi_1\circ (u')^a\quad\text{weakly}^{\ast}\text{ in }BV({]a,b[})\,,\label{wPsi}\\
& u_{k}^{\prime}\rightarrow(u^{\prime})^{a}\quad\text{pointwise }\mathcal{L}%
^{1}\text{-a.e.\ in }{]a,b[}\,.\label{conv l1}%
\end{align}
\end{theorem}

\begin{proof}
By  (\ref{functional}) and (\ref{supenergy}) we have  that each $u_{k}$
belongs to $W^{2,1}({]a,b[})$ and
\begin{equation}
C_{1}:=\sup_{k}\int_{a}^{b}[\,|u_{k}|+|u_{k}^{\prime}|+\psi(u_{k}^{\prime
})|u_{k}^{\prime\prime}|\,]\,dx<+\infty\,.\label{BoundC}%
\end{equation}
Let us define
\begin{equation}
v_{k}:=\Psi_1\circ u_{k}^{\prime}\,.\label{vk}%
\end{equation}
As $\Psi_1$ is Lipschitz in ${\mathbb{R}}$, the functions $v_{k}$
belong to $W^{1,1}({]a,b[})$ and
\begin{equation}
v_{k}^{\prime}=\psi(u_{k}^{\prime})u_{k}^{\prime\prime}\qquad{\mathcal{L}%
}^{1}\text{-a.e.\ on }{]a,b[}\,.\label{derivvk}%
\end{equation}
It follows from  (\ref{intfin}) and (\ref{supenergy}) that
\begin{equation}
\int_{a}^{b}[\,|v_{k}|+|v_{k}^{\prime}|\,]\,dx\leq M(b-a)+C\,.\label{BoundC2}%
\end{equation}
By (\ref{BoundC}) and (\ref{BoundC2}) and the Helly theorem, passing
to a subsequence if necessary, we may assume that
$$
u_{k}\rightharpoonup u\quad\text{weakly}^{\ast}\text{ in
}BV({]a,b[})
$$
and
\begin{equation}
v_{k}\left(  x\right)  \rightarrow v\left(  x\right)  \quad\text{for all }%
x\in{]a,b[}\label{convpoint3}%
\end{equation}
for some  $u\in BV({]a,b[})$ and  $v\colon
{]a,b[}\rightarrow[0,M]$ with pointwise bounded variation. Note that
(\ref{convpoint3}) determines the values of $v$ at every $x\in{]a,b[}$.

Since $\Psi^{-1}_1$ is continuous, we obtain
\begin{equation}
u_{k}^{\prime}\rightarrow w:=\Psi^{-1}_1(v)\quad\text{pointwise in
}{]a,b[}\,.
\label{convpoint4}%
\end{equation}
Moreover $w$ has left and right limits in $\overline
{\mathbb{R}}$ at each point $x\in{]a,b[}$, denoted by $w_-(x)$ and $w_+(x)$
respectively, and
\begin{equation}
w(x)=w_-(x)=w_+(x)\quad\text{except for a countable set of points }x\,.
\label{countexcept}%
\end{equation}
We now split the remaining part of the proof into two steps.

\noindent\textbf{Step 1: }We prove that
\begin{equation}
w=(u^{\prime})^{a}\quad{\mathcal{L}}^{1}\text{-a.e.\ in }{]a,b[}\,.
\label{w=dotu}%
\end{equation}
If not, we have ${\mathcal{L}}^{1}(\{w\neq(u^{\prime})^{a}\})>0$. By
(\ref{inf>0}) the function $\Psi^{-1}_1$ is locally Lipschitz and so
$w=\Psi^{-1}_1(v)$ is finite ${\mathcal{L}}^{1}$-a.e.\ since $v\in L^{1}%
({]a,b[})$. Hence there exists $t_{0}>0$ such that
\[
{\mathcal{L}}^{1}(\{w\neq(u^{\prime})^{a}\}\cap\{|w|<t_{0}\})>0\,
\]
and, in particular, we may find an infinite number of disjoint open intervals
$I$ such that
\begin{equation}
{\mathcal{L}}^{1}(\{w\neq(u^{\prime})^{a}\}\cap\{|w|<t_{0}\}\cap I)>0\,.
\label{wneqdotu}%
\end{equation}
By a change of variables we obtain
\begin{equation}
\int_{I}\psi(u_{k}^{\prime})|u_{k}^{\prime\prime}|\,dx\geq\int_{m_{k}}^{M_{k}%
}\psi(t)\,dt\,, \label{changevar2}%
\end{equation}
where
\[
m_{k}:=\inf_{I}u_{k}^{\prime}\quad\text{and}\quad M_{k}:=\sup_{I}u_{k}%
^{\prime}\,.
\]
We claim that at least one of the two sequences $\{m_{k}\}$ and
$\{M_{k}\}$ is divergent. Indeed, if not, a subsequence of
$\{u_{k}^{\prime}\}$ would be bounded in $L^{\infty}(I)$. This
implies that $u^{\prime}\in L^{\infty}(I)$ and that
$u_{k}^{\prime}\rightharpoonup u^{\prime}$ weakly$^{\ast}$ in
$L^{\infty}(I)$. As $u_{k}^{\prime}\rightarrow w$ pointwise
${\mathcal{L}}^{1}$-a.e.\ in $I$, we deduce that $u^{\prime}=w$
${\mathcal{L}}^{1}$-a.e.\ in $I$, which
contradicts~\eqref{wneqdotu}. Hence the claim holds. If%
\begin{equation}
\lim_{k\rightarrow\infty}\,M_{k}=+\infty\,, \label{sup2}%
\end{equation}
then by \eqref{convpoint4} and \eqref{wneqdotu}
\begin{equation}
\limsup_{k\rightarrow\infty}\,m_{k}<t_{0}\,. \label{inf2}%
\end{equation}
From \eqref{changevar2}, \eqref{inf2}, and \eqref{sup2}  we obtain
\[
\liminf_{k\rightarrow\infty}\int_{I}\psi(u_{k}^{\prime})|u_{k}^{\prime\prime
}|\,dx\geq\int_{t_{0}}^{+\infty}\psi(t)\,dt>0\,.
\]
Analogously, if $\lim_{k}\,m_{k}=-\infty$ then
\[
\liminf_{k\rightarrow\infty}\int_{I}\psi(u_{k}^{\prime})|u_{k}^{\prime\prime
}|\,dx\geq\int_{-\infty}^{-t_{0}}\psi(t)\,dt>0\,.
\]
In any case we can choose an arbitrarily large number $m$ of
disjoint intervals $I$ satisfying \eqref{wneqdotu}. Adding the
contributions of each interval we obtain
\[
\liminf_{k\rightarrow\infty}\int_{a}^{b}\psi(u_{k}^{\prime})|u_{k}%
^{\prime\prime}|\,dx\geq m\min\Big\{  \int_{t_{0}}^{+\infty}\psi
(t)\,dx,\int_{-\infty}^{-t_{0}}\psi(t)\,dt\Big\}  \,,
\]
which contradicts~\eqref{BoundC}  for $m$ large enough. This concludes the
proof of~\eqref{w=dotu}  .

\noindent\textbf{Step 2:} To prove that $u\in X^1_{\psi}({]a,b[})$ it remains to
show that the positive part $\left(  (u^{\prime})^{c}\right)  ^{+}$ and the
negative part $\left(  (u^{\prime})^{c}\right)  ^{-}$ of the measure
$(u^{\prime})^{c}$ are concentrated on $Z^{+}[(u^{\prime})^{a}]$ and
$Z^{-}[(u^{\prime})^{a}]$ respectively, that is%
\begin{equation}
\left(  (u^{\prime})^{c}\right)  ^{\pm}({]a,b[}\setminus Z^{\pm}[(u^{\prime
})^{a}])=0\,.\label{cantor}%
\end{equation}
To this purpose we introduce the sets
\begin{align}
E^{+}[u^{\prime}] &  :=\Big\{ x\in{]a,b[}:\lim_{\varepsilon\rightarrow
0+}\frac{\left(  u^{\prime}\right)  ^{+}({]x-\varepsilon,x+\varepsilon\lbrack
})}{2\varepsilon}=+\infty\Big\} \,,\label{Eu+}\\
E^{-}[u^{\prime}] &  :=\Big\{ x\in{]a,b[}:\lim_{\varepsilon\rightarrow
0+}\frac{\left(  u^{\prime}\right)  ^{-}({]x-\varepsilon,x+\varepsilon\lbrack
})}{2\varepsilon}=+\infty\Big\} \,,\label{Eu-}\\
E[u^{\prime}] &  :=\Big\{ x\in{]a,b[}:\lim_{\varepsilon\rightarrow
0+}\frac{|u^{\prime}|({]x-\varepsilon,x+\varepsilon\lbrack})}{2\varepsilon
}=+\infty\Big\} \,.\nonumber%
\end{align}
Since $((u^{\prime})^{s})^+=((u^{\prime})^{+})^s$ is concentrated on $E^{+}[u^{\prime}]$ and
$((u^{\prime})^{s})^-=((u^{\prime})^{-})^s$ is concentrated on $E^{-}[u^{\prime}]$
(see, e.g., \cite[Theorem 2.22]{AFP}), to prove \eqref{cantor}  it is enough
to show that
\begin{equation}
E^{+}[u^{\prime}]\setminus Z^{+}[(u^{\prime})^{a}]\text{ and }E^{-}[u^{\prime
}]\setminus Z^{-}[(u^{\prime})^{a}]\text{ are at most countable.}%
\label{countable}%
\end{equation}
We only show that $E^{+}[u^{\prime}]\setminus
Z^{+}[(u^{\prime})^{a}]$ is at most countable, since the other
property can be proved in a similar way. Assume by contradiction
that $E^{+}[u^{\prime}]\setminus Z^{+}[(u^{\prime
})^{a}]$ is not countable. Since by \eqref{Zu} and \eqref{w=dotu}%
\[
Z^{+}[(u^{\prime})^{a}]\subset\{x\in{]a,b[}:\max\{w_-(x),w_+(x)\}=+\infty\}\,,
\]
by \eqref{countexcept}  there exists $t_{0}>0$ such that
\[
(E^{+}[u^{\prime}]\setminus Z^{+}[(u^{\prime})^{a}])\cap\{w<t_{0}%
\}\quad\text{is uncountable.}%
\]
 Fix $t_{1}>t_{0}$ and let
 $x_{1},\dots,x_{m}$ be $m$ distinct points in $(E^{+}[u^{\prime
}]\setminus Z^{+}[(u^{\prime})^{a}])\cap\{w<t_{0}\}$. By \eqref{Eu+}  there
exists $\varepsilon>0$ such that the intervals $I_{j}:={]x_{j}-\varepsilon
,x_{j}+\varepsilon\lbrack}$ are pairwise disjoint and
\begin{equation}
\frac{\left(  u^{\prime}\right)  ^{+}({]x_{j}-\varepsilon,x_{j}+\varepsilon
\lbrack})}{2\varepsilon}>t_{1}\quad\text{for }i=1,\dots,m\,. \label{quotient}%
\end{equation}
By a change of variables we obtain
\begin{equation}
\int_{I_{j}}\psi(u_{k}^{\prime})|u_{k}^{\prime\prime}|\,dx\geq\int_{m_{kj}%
}^{M_{kj}}\psi(t)\,dt\,, \label{changevar10}%
\end{equation}
where
\[
m_{kj}:=\inf_{I_{j}}u_{k}^{\prime}\quad\text{and}\quad M_{kj}:=\sup_{I_{j}%
}u_{k}^{\prime}\,.
\]
By \eqref{convpoint4}  and the fact that $w\left(  {x_{j}}\right)  <t_{0}$ we
deduce that%
\begin{equation}
\limsup_{k\rightarrow\infty}\,m_{kj}<t_{0} \label{inf20}%
\end{equation}
for $j=1,\dots,m$. On the other hand,  \eqref{convw*}  and
\eqref{quotient} yield
$$
\liminf_{k\rightarrow\infty}\frac{1}{2\varepsilon}\int_{x_{j}-\varepsilon
}^{x_{j}+\varepsilon}(u_{k}^{\prime})^{+}\,dx  \geq \frac{\left(  u^{\prime}\right)  ^{+}({]x_{j}-\varepsilon,x_{j}%
+\varepsilon\lbrack})}{2\varepsilon}>t_{1}
$$
(this can be seen as a particular case of the Reshetnyak lower semicontinuity theorem, with $f=\left(  \cdot\right)
^{+}$). This  implies that
$\liminf_{k\rightarrow\infty}\,M_{kj}>t_{1}$
for $j=1,\dots,m$. Hence, also by  \eqref{changevar10}   and
\eqref{inf20},    we obtain
\[
\liminf_{k\rightarrow\infty}\sum_{j=1}^{m}\,\int_{I_{j}}\psi(u_{k}^{\prime
})|u_{k}^{\prime\prime}|\,dx\geq\sum_{j=1}^{m}\,\liminf_{k\rightarrow\infty
}\int_{I_{j}}\psi(u_{k}^{\prime})|u_{k}^{\prime\prime}|\,dx\geq m\int_{t_{0}%
}^{t_{1}}\psi(t)\,dt\,,
\]
which contradicts~\eqref{BoundC} for $m$ large enough. This shows (\ref{countable}) and concludes
the proof of the theorem.
\end{proof}

\subsection{Relaxation}

The following theorem, which is the main result of the section, is devoted to the characterization of the
relaxation of $\mathcal{F}_1$ with respect to strong
convergence in $L^{1}({]a,b[})$.

\begin{theorem}
\label{th:relaxation} Let
$\overline{\mathcal{F}}_1\colon L^{1}({]a,b[}%
)\rightarrow\lbrack0,+\infty]$ be defined by:
\begin{equation*}
\overline{\mathcal{F}}_1\left(  u\right)  :=\inf\left\{  \liminf_{k\rightarrow
\infty}\mathcal{F}_1\left(  u_{k}\right)  :\,u_{k}\rightarrow u\text{ in }%
L^{1}({]a,b[})\right\}%
\end{equation*}
for every $u\in L^{1}({]a,b[})$.
Then
\begin{equation}
\overline{\mathcal{F}}_1(u)=%
\begin{cases}
\displaystyle |u^{\prime}|({]a,b[})+|v^{\prime}|({]a,b[}\setminus S_{u}%
)+\sum_{x\in S_{u}}\Phi(\nu_{u},(u^{\prime})_{-}^{a},(u^{\prime})_{+}^{a})\, &
\text{if }u\in X^1_{\psi}({]a,b[})\,,\\
+\infty & \text{otherwise,}%
\end{cases}
\label{relaxed}%
\end{equation}
where $v:=\Psi_1\circ (u')^a$, $\nu_{u}:=\mathrm{sign}(u_{+}-u_{-})$, and
\begin{equation}%
\begin{array}
[c]{c}%
\displaystyle \Phi(1,t_{1},t_{2}):=\int_{t_{1}}^{+\infty}\psi(t)\,dt+\int
_{t_{2}}^{+\infty}\psi(t)\,dt\,,\\
\displaystyle \Phi(-1,t_{1},t_{2}):=\int_{-\infty}^{t_{1}}\psi(t)\,dt+\int
_{-\infty}^{t_{2}}\psi(t)\,dt\,.
\end{array}
\label{Phi}%
\end{equation}
\end{theorem}

\begin{remark}
\label{rem2.5'}
For every $x\in S_{u}$ we have
$$
\Phi(\nu_{u}(x),(u^{\prime})_{-}^{a}(x),(u^{\prime})_{+}^{a}(x))
=|v'|(\{x\})+\hat \Phi(\nu_{u}(x),(u^{\prime})_{-}^{a}(x),(u^{\prime})_{+}^{a}(x))\,,
$$
where
$$
\displaystyle \hat\Phi(1,t_{1},t_{2}):=\int_{\max\{t_{1},t_{2}\}}^{+\infty}\psi(t)\,dt\qquad
\text{and}\qquad
\displaystyle \hat\Phi(-1,t_{1},t_{2}):=\int_{-\infty}^{\min\{t_{1},t_{2}\}}\psi(t)\,dt\,.
$$
In particular, for every Borel set $B\subset{]a,b[}$
$$
|v^{\prime}|(B\setminus S_{u}%
)+\sum_{x\in S_{u}\cap B}\Phi(\nu_{u},(u^{\prime})_{-}^{a},(u^{\prime})_{+}^{a})
= |v^{\prime}|(B%
)+\sum_{x\in S_{u}\cap B}\hat\Phi(\nu_{u},(u^{\prime})_{-}^{a},(u^{\prime})_{+}^{a})\ge
|v^{\prime}|(B)\,.
$$
\end{remark}

\begin{proof}[Proof of Theorem~\ref{th:relaxation}]
Let $\mathcal{G}$ be the functional defined by the right hand side of
(\ref{relaxed}). We prove that for every $u_{k}\rightarrow u$ in
$L^{1}({]a,b[})$ we have
\begin{equation}
\mathcal{G}(u)\leq\liminf_{k\rightarrow\infty}{\mathcal{F}}_1(u_{k}%
)\,.\label{liminf}%
\end{equation}
It is enough to consider sequences $\{u_{k}\}$ for which the liminf
is a limit and has a finite value and $u_k\to u$ pointwise
$\Lone$-a.e.\ in ${]a,b[}$. Then $u_{k}$ belongs to
$W^{2,1}({]a,b[})$ and \eqref{supenergy}  is satisfied. This implies
that
\begin{equation}
|u^{\prime}|({]a,b[})\leq\liminf_{k\rightarrow\infty}\int_{a}^{b}
|u_{k}^{\prime}|\,dx\,.\label{lscu'}
\end{equation}
Moreover, it follows from Theorem \ref{theorem compactness} that
$u\in X_{\psi}^1({]a,b[})$ and that, up to a subsequence,
$\{u_{k}^{\prime }\}$ converges to $(u^{\prime})^{a}$ pointwise
${\mathcal{L}}^{1}$-a.e.\ in ${]a,b[}$.

Let $F$ be a finite subset of $S_{u}$. We want to prove that
\begin{equation}
|v^{\prime}|({]a,b[}\setminus F)+\sum_{x\in F}\Phi(\nu_{u},(u^{\prime}
)_{-}^{a},(u^{\prime})_{+}^{a})\leq\liminf_{k\rightarrow\infty}\int_{a}
^{b}\psi(u_{k}^{\prime})|u_{k}^{\prime\prime}|\,dx\,. \label{lscv'}
\end{equation}
We write $F$ as $\{x_{1},\dots,x_{m}\}$, with $a<x_{1}<\dots<x_{m}<b$. For
every $\varepsilon>0$ there exists $\delta=\delta(\varepsilon)\in{]0,\varepsilon\lbrack}$ such
that $a<x_{1}-\delta<x_{1}+\delta<x_{2}-\delta<x_{2}+\delta<\dots
<x_{m-1}-\delta<x_{m-1}+\delta<x_{m}-\delta<x_{m}+\delta<b$ and
\begin{align}
&  |u(x_{j}-\delta)-u_{-}(x_{j})|<\varepsilon\,,\quad|u(x_{j}+\delta
)-u_{+}(x_{j})|<\varepsilon\,,\label{u-delta}\\
&  |(u^{\prime})^{a}(x_{j}-\delta)-(u^{\prime})_{-}^{a}(x_{j})|<\varepsilon
\,,\quad|(u^{\prime})^{a}(x_{j}+\delta)-(u^{\prime})_{+}^{a}(x_{j}%
)|<\varepsilon\,,\label{uprime-delta}\\
&  u_{k}(x_{j}-\delta)\rightarrow u(x_{j}-\delta)\,,\quad u_{k}(x_{j}%
+\delta)\rightarrow u(x_{j}+\delta)\quad\text{as}\quad k\rightarrow
\infty\,,\label{convu-delta}\\
&  u_{k}^{\prime}(x_{j}-\delta)\rightarrow(u^{\prime})^{a}(x_{j}%
-\delta)\,,\quad u_{k}^{\prime}(x_{j}+\delta)\rightarrow(u^{\prime})^{a}%
(x_{j}+\delta)\quad\text{as}\quad k\rightarrow\infty\,,
\label{convuprime-delta}%
\\
&  |(u^{\prime})^{a}(x_{j}-\delta)| +
|(u^{\prime})^{a}(x_{j}+\delta)|+\varepsilon <
\frac{|[u](x_j)|-4\varepsilon}{2\delta}\,,\nonumber
\end{align}
for $j=1,\dots,m$.

Since $v_{k}\rightarrow v$ pointwise ${\mathcal{L}}^{1}$-a.e.\ in ${]a,b[}$
and $v_{k}^{\prime}=\psi(u_{k}^{\prime})u_{k}^{\prime\prime}$ ${\mathcal{L}%
}^{1}$-a.e.\ in ${]a,b[}$, we obtain
\[
|v^{\prime}|({]x_{j}+\delta,x_{j+1}-\delta\lbrack})\leq\liminf_{k\rightarrow
\infty}\int_{x_{j}+\delta}^{x_{j+1}-\delta}\psi(u_{k}^{\prime})|u_{k}%
^{\prime\prime}|\,dx
\]
for $j=1,\dots,m-1$. A similar result holds for the intervals ${]a,x_{1}%
-\delta\lbrack}$ and ${]x_{m}+\delta,b[}$. Let $F_{\delta}$ be the union of
the intervals $[x_{j}-\delta,x_{j}+\delta]$ for $j=1,\dots,m$. Summing with
respect to $j$, and adding the contributions of the intervals ${]a,x_{1}%
-\delta\lbrack}$ and ${]x_{m}+\delta,b[}$, we obtain
\begin{equation}
|v^{\prime}|({]a,b[}\setminus F_{\delta})\leq\liminf_{k\rightarrow\infty}%
\int_{{]a,b[}\setminus F_{\delta}}\psi(u_{k}^{\prime})|u_{k}^{\prime\prime
}|\,dx\,. \label{vprime}%
\end{equation}

We consider now the interval $I_{j}^{\delta}:=[x_{j}-\delta,x_{j}+\delta]$,
assuming that $[u](x_{j})=u_{+}(x_{j})-u_{-}(x_{j})>0$. By the mean value
theorem there exists $y_{kj}^{\delta}\in{]x_{j}-\delta,x_{j}+\delta\lbrack}$
such that
\begin{equation}
u_{k}^{\prime}(y_{kj}^{\delta})=\frac{u_{k}(x_{j}+\delta)-u_{k}(x_{j}-\delta
)}{2\delta}\geq\frac{\lbrack u](x_{j})-4\varepsilon}{2\delta}\,,
\label{uprimekyj}%
\end{equation}
where the last inequality follows from \eqref{u-delta}  and~(\ref{convu-delta}%
) for $k$ sufficiently large. By a change of variables we obtain
\begin{align*}
&  \displaystyle  \int_{u_{k}^{\prime}(x_{j}-\delta)}^{\frac{[u](x_{j}%
)-4\varepsilon}{2\delta}}\psi(t)\,dt\leq\int_{x_{j}-\delta}^{y_{kj}^{\delta}%
}\psi(u_{k}^{\prime})|u_{k}^{\prime\prime}|\,dx\,,\\
&  \displaystyle  \int_{u_{k}^{\prime}(x_{j}+\delta)}^{\frac{[u](x_{j}%
)-4\varepsilon}{2\delta}}\psi(t)\,dt\leq\int_{y_{kj}^{\delta}}^{x_{j}+\delta
}\psi(u_{k}^{\prime})|u_{k}^{\prime\prime}|\,dx\,.
\end{align*}
Adding these inequalities and taking the limit as $k\rightarrow\infty$ we
obtain, thanks to~\eqref{convuprime-delta}  ,
\begin{equation}
\int_{(u^{\prime})^{a}(x_{j}-\delta)}^{\frac{[u](x_{j})-4\varepsilon}{2\delta
}}\psi(t)\,dt+\int_{(u^{\prime})^{a}(x_{j}+\delta)}^{\frac{[u](x_{j}%
)-4\varepsilon}{2\delta}}\psi(t)\,dt\leq\liminf_{k\rightarrow\infty}%
\int_{x_{j}-\delta}^{x_{j}-\delta}\psi(u_{k}^{\prime})|u_{k}^{\prime\prime
}|\,dx\,. \label{jumps-delta+}%
\end{equation}
Similarly, if $[u](x_{j})<0$, then we have
\begin{equation}
\int_{\frac{\lbrack u](x_{j})+4\varepsilon}{2\delta}}^{(u^{\prime})^{a}%
(x_{j}-\delta)}\psi(t)\,dt+\int_{\frac{[u](x_{j})+4\varepsilon}{2\delta}%
}^{(u^{\prime})^{a}(x_{j}+\delta)}\psi(t)\,dt\leq\liminf_{k\rightarrow\infty
}\int_{x_{j}-\delta}^{x_{j}-\delta}\psi(u_{k}^{\prime})|u_{k}^{\prime\prime
}|\,dx\,. \label{jumps-delta-}%
\end{equation}
From \eqref{vprime}, \eqref{jumps-delta+}, and \eqref{jumps-delta-}
we deduce that
\begin{equation*}%
\begin{array}
[c]{c}%
\displaystyle  |v^{\prime}|({]a,b[}\setminus F_{\delta})+\sum_{[u](x_{j}%
)>0}\Big(  \int_{(u^{\prime})^{a}(x_{j}-\delta)}^{\frac{[u](x_{j}%
)-4\varepsilon}{2\delta}}\psi(t)\,dt+\int_{(u^{\prime})^{a}(x_{j}+\delta
)}^{\frac{[u](x_{j})-4\varepsilon}{2\delta}}\psi(t)\,dt\Big)\\
\displaystyle  {}+\sum_{[u](x_{j})<0}\Big(  \int_{\frac{[u](x_{j}%
)+4\varepsilon}{2\delta}}^{(u^{\prime})^{a}(x_{j}-\delta)}\psi(t)\,dt+\int
_{\frac{[u](x_{j})+4\varepsilon}{2\delta}}^{(u^{\prime})^{a}(x_{j}+\delta
)}\psi(t)\,dt\Big)\\
\displaystyle  \leq\liminf_{k\rightarrow\infty}
\int_{a}^{b}\psi
(u_{k}^{\prime})|u_{k}^{\prime\prime}|\,dx\,.
\end{array}
\end{equation*}
Taking the limit as $\varepsilon\rightarrow0$ (which implies
$\delta(\varepsilon)\rightarrow0$)
we obtain \eqref{lscv'}  thanks to~\eqref{uprime-delta}  .

Since $S_{u}$ is at most countable, \eqref{liminf}  can be obtained from
\eqref{lscv'}  by taking the supremum over all finite sets $F$ contained in
$S_{u}$.

Conversely, let $u\in X^1_{\psi}({]a,b[})$. We claim that there
exists a sequence $\{u_k\}$ in $W^{2,1}\left(  \left]  a,b\right[
\right)  $ such that $u_{k}\rightarrow
u$ in $L^{1}\left(  \left]  a,b\right[  \right)  $ and%
\begin{equation}
\mathcal{G}\left(  u\right)  \geq\limsup_{k\rightarrow\infty}\mathcal{F}_1
\left(  u_{k}\right) \, .\label{limsup inequality}%
\end{equation}
It is clearly enough to consider the case
$\mathcal{G}\left(  u\right)  <+\infty$.

We divide the proof into three steps.

\noindent\textbf{Step 1:} We prove (\ref{limsup inequality}) under the
additional assumptions that $\left(  u^{\prime}\right)  ^{a}$ is bounded and
that  $S_{u}=\left\{  x_{1},\ldots,x_{m}\right\}  $, with
$x_{1}<\ldots<x_{m}$. Note that in this case $Z[(u^{\prime})^{a}]=\emptyset$, hence
$(u^{\prime})^{c}=0$.

Construct a sequence $\{v_k\}$ in $W^{1,1}\left(  \left]  a,b\right[
\right)  $ such that $v_{k}\rightarrow v=\Psi_1\circ (  u^{\prime})
^{a}$
 pointwise $\Lone$-a.e.\ in $\left]  a,b\right[  $, $\Psi_1\left(  -\left\|  \left(
u^{\prime}\right)  ^{a}\right\|  _{\infty}\right)  \leq v_{k}\leq\Psi_1\left(
\left\|  \left(  u^{\prime}\right)  ^{a}\right\|  _{\infty}\right)  $, and%
\begin{equation*}
\int_{a}^{b}\left|  v_{k}^{\prime}\left(  x\right)  \right|
\,dx\rightarrow \left|  v^{\prime}\right|  \left(  \left] a,b\right[
\right)\,  .
\end{equation*}
Setting $w_{k}:=\Psi_1^{-1}\left(  v_{k}\right)  $, we have $w_{k}\in
W^{1,1}({]a,b[})$ thanks to \eqref{inf>0},
\begin{equation}\label{5689}
w_{k}\rightarrow\left(
u^{\prime}\right)  ^{a}\quad  \text{pointwise }{\mathcal{L}}^{1}\text{-a.e.\ in }
{]a,b[}\,,
\end{equation}
and $\left\|  w_{k}\right\|  _{\infty}\leq\left\|  \left(
u^{\prime}\right)  ^{a}\right\|  _{\infty}$. Find $\delta_{k}\rightarrow0^{+}$
such that
\begin{equation}
w_{k}\left(  x_{j}-\delta_{k}\right)  \rightarrow\left(  u^{\prime}\right)
_{-}^{a}\left(  x_{j}\right)  ,\quad w_{k}\left(  x_{j}+\delta_{k}\right)
\rightarrow\left(  u^{\prime}\right)  _{+}^{a}\left(  x_{j}\right)\quad\text{for $j=1,\ldots,m$,}
\label{delta k}%
\end{equation}
 and
 \begin{equation}\label{deltak2}
 \begin{array}{c}
 \displaystyle\int_{x_{j-1}+\delta_k}^{x_j-\delta_k}|v'_k|\, dx\to |v'|(]x_{j-1},x_j[)
 \quad \text{for $j=2,\ldots,m$,}\smallskip\\
 \displaystyle \int_{a}^{x_1-\delta_k}|v'_k|\, dx\to |v'|(]a,x_1[)\,,\  \int_{x_m+\delta_k}^{b}|v'_k|\, dx\to |v'|(]x_m+\delta_k,b[)\,.
 \end{array}
 \end{equation}
By \eqref{5689} and by the dominated convergence theorem we have
\begin{equation}\label{6689}
u_+(x_{j-1})+\int_{x_{j-1}+\delta_k}^{x_j-\delta_k} w_{k}(s)\, ds \longrightarrow
u_+(x_{j-1})+\int_{x_{j-1}}^{x_j} (u^{\prime})^{a}\, ds
= u_-(x_j)
\end{equation}
for $j=2,\ldots,m$, with the obvious changes for $j=1$ and $j={m+1}$.

To deal with the jump point $x_j$, assume first that
\begin{equation}\label{case1}
u_{+}\left(  x_{j}\right)-u_{-}\left(  x_{j}\right)>0\,.
\end{equation}
In this case we need to construct functions
$f_{kj}\in C^{2}(  \left[
x_{j}-\delta_{k},x_{j}+\delta_{k}\right])$ that satisfy the following properties:
there exist $y_{kj}\in{]x_j-\delta_k,x_j+\delta_k[}$ such that
\begin{eqnarray}
&\displaystyle f_{kj}\left(  x_{j}-\delta_{k}\right)     =u_{+}\left(  x_{j-1}\right)
+\int_{x_{j-1}+\delta_{k}}^{x_{j}-\delta_{k}}w_{k}\left(  s\right)  \,ds,\quad
f_{kj}\left(  x_{j}+\delta_{k}\right)  =u_{+}\left(  x_{j}\right)
,\label{initial condition I}\\
& f_{kj}^{\prime}\left(  x_{j}-\delta_{k}\right)  =w_{k}\left(  x_{j}
-\delta_{k}\right)  ,\quad f_{kj}^{\prime}\left(  x_{j}+\delta_{k}\right)
=w_{k}\left(  x_{j}+\delta_{k}\right)\,,\label{initial condition II}\\
&\displaystyle\vphantom{\int}\!\!\!\!\!\!\!\! f''_{kj}(x)>0\text{ for  }x\in{]x_j-\delta_k,y_{kj}[}\,,\quad
 f''_{kj}(x)<0\text{ for }x\in{]y_{kj}, x_j+\delta_k[}\,,
\label{1inflection}\\
&   \displaystyle |f_{kj}\left(  x_{j}-\delta_{k}\right)-\min_{[
x_{j}-\delta_{k},y_{kj}]}f_{kj}|\leq\tfrac 1k\,,\quad
|f_{kj}\left(  x_{j}+\delta_{k}\right)-\max_{[
y_{kj},x_{j}+\delta_{k}]}f_{kj}|\leq\tfrac 1k\,, \label{minmax}
\end{eqnarray}
where we replace $x_{j-1}$ and $x_{j-1}-\delta_{k}$ by $a$ in the case $j=1$.

We now discuss briefly the existence of such functions. We observe that the latter conditions in equations \eqref{initial condition I}--\eqref{1inflection} imply that the
graph of
$f_{kj}$  in the interval ${[y_{kj},x_j+\delta_k[}$ lies below the straight line passing through the point $(x_j+\delta_k,u_+(x_j))$ with slope $w_k(x_j+\delta_k)$, i.e.,
$$
f_{kj}(x)\le u_+(x_j)+w_k(x_j+\delta_k)(x-x_j-\delta_k)
$$
for $x\in {[y_{kj},x_j+\delta_k[}$. It is then easy to see that the inequality
\begin{equation}\label{4356}
u_{+}\left(  x_{j}\right)-2w_{k}\left(  x_{j}+\delta_{k}\right)\delta_k-u_{+}\left(  x_{j-1}\right)
-\int_{x_{j-1}+\delta_{k}}^{x_{j}-\delta_{k}}w_{k}\left(  s\right)  \,ds>0\,,
\end{equation}
allows to fulfill also the former conditions in equations
\eqref{initial condition I}--\eqref{1inflection}, as well as \eqref{minmax}.
By \eqref{delta k}, \eqref{6689}, and
\eqref{case1}, inequality
\eqref{4356} is satisfied when $\delta_k$ is small enough.

If the left-hand side of \eqref{case1} is negative then we choose $f_{kj}$ so that
\eqref{initial condition I}
and \eqref{initial condition II} hold, and there exists  $y_{kj}\in {]x_j-\delta_k, x_j+\delta_k[}$ such that
\begin{eqnarray*}
& \!\!\!\!\!\!\!\! f''_{kj}(x)<0\text{ for  }x\in{]x_j-\delta_k,y_{kj}[}\,,\quad
 f''_{kj}(x)>0\text{ for }x\in{]y_{kj}, x_j+\delta_k[}\,, \\
& \displaystyle  |f_{kj}\left(  x_{j}-\delta_{k}\right)-\max_{[
x_{j}-\delta_{k},x_{j}+\delta_{k}]}f_{kj}|\leq\tfrac 1k\,,\quad
|f_{kj}\left(  x_{j}+\delta_{k}\right)-\min_{[
x_{j}-\delta_{k},x_{j}+\delta_{k}]}f_{kj}|\leq\tfrac 1k\,.
\end{eqnarray*}
In the same way the construction is possible if $\delta_k$ is small enough.

We are now ready to define the approximating sequence
\[
u_{k}\left(  x\right)  :=\left\{
\begin{array}
[c]{ll}%
u_{+}\left(  a\right)  +\int_{a}^{x}w_{k}\left(  s\right)  \,ds & \text{if
}a\leq x<x_{1}-\delta_{k}\,,\\
f_{kj}\left(  x\right)  & \text{if }x_{j}-\delta_{k}\leq x<x_{j}+\delta
_{k},\text{ }j=1,\ldots,m\,,\\
u_{+}\left(  x_{j-1}\right)  +\int_{x_{j-1}+\delta_{k}}^{x}w_{k}\left(
s\right)  \,ds & \text{if }x_{j-1}+\delta_{k}\leq x<x_{j}-\delta_{k},\text{
}j=2,\ldots,m\,,\\
u_{+}\left(  x_{m}\right)  +\int_{x_{m}+\delta_{k}}^{x}w_{k}\left(
s\right) \,ds & \text{if }x_{m}+\delta_{k}\leq x<b\,.
\end{array}
\right.
\]
Let us define $x_0:=a$ and $x_{m+1}:=b$. Since $w_{k}\rightarrow\left(  u^{\prime}\right)  ^{a}$ in $L^{1}\left(
\left]  a,b\right[  \right)  $, we have
\[
u_{k}\left(  x\right)  \rightarrow u_{+}\left(  x_{j-1}\right)  +\int
_{x_{j-1}}^{x}\left(  u^{\prime}\right)  ^{a}\left(  s\right)  \,ds=u\left(
x\right)
\]
for every $x\in{]x_{j-1},x_{j}[}$ and $j=1,\ldots,m+1$
and, in turn,  $u_{k}\rightarrow u$ in $L^{1}\left(  \left]
a,b\right[  \right)  $. As%
\begin{align*}
\int_{x_{j-1}+\delta_{k}}^{x_{j}-\delta_{k}}\left|  u_{k}^{\prime}\right|
\,dx+\int_{x_{j-1}+\delta_{k}}^{x_{j}-\delta_{k}}\psi\left(  u_{k}^{\prime
}\right)  \left|  u_{k}^{\prime\prime}\right|  \,dx&=\int_{x_{j-1}+\delta_{k}%
}^{x_{j}-\delta_{k}}\left|  w_{k}\right|  \,dx+\int_{x_{j-1}+\delta_{k}%
}^{x_{j}-\delta_{k}}\psi\left(  w_{k}\right)  \left|  w_{k}^{\prime}\right|
\,dx\\
&\leq\int_{x_{j-1}}^{x_{j}}\left|  w_{k}\right|
\,dx+\int_{x_{j-1}+\delta_k}^{x_{j}-\delta_k }\left|
v_{k}^{\prime}\right|  \,dx\,,
\end{align*}
by (\ref{deltak2}) and the fact that $w_{k}\rightarrow\left(  u^{\prime
}\right)  ^{a}$ in $L^{1}\left(  \left]  a,b\right[  \right)  $ we have%
\begin{equation}
\begin{array}{c}
\displaystyle \limsup_{k\rightarrow\infty}\Big(  \int_{x_{j-1}+\delta_{k}}^{x_{j}
-\delta_{k}}\left|  u_{k}^{\prime}\right|  \,dx+\int_{x_{j-1}+\delta_{k}
}^{x_{j}-\delta_{k}}\psi\left(  u_{k}^{\prime}\right)  \left|  u_{k}
^{\prime\prime}\right|  \,dx\Big)\smallskip\\
\displaystyle\leq \int_{x_{j-1}}^{x_{j}}\left|  \left(
u^{\prime}\right)  ^{a}\right|  \,dx+\left|  v^{\prime}\right|  \left(
\left]  x_{j-1},x_{j}\right[  \right)\,. \label{lim1}
\end{array}
\end{equation}
Similarly,
\begin{align}
\limsup_{k\rightarrow\infty}\Big(  \int_{a}^{x_{1}-\delta_{k}}\left|
u_{k}^{\prime}\right|  \,dx+\int_{a}^{x_{1}-\delta_{k}}\psi\left(
u_{k}^{\prime}\right)  \left|  u_{k}^{\prime\prime}\right|  \,dx\Big)   &
\leq\int_{a}^{x_{1}}\left|  \left(  u^{\prime}\right)  ^{a}\right|
\,dx+\left|  v^{\prime}\right|  \left(  \left]  a,x_{1}\right[  \right)
\,,\label{lim2}\\
\limsup_{k\rightarrow\infty}\Big(  \int_{x_{m}+\delta_{k}}^{b}\left|
u_{k}^{\prime}\right|  \,dx+\int_{x_{m}+\delta_{k}}^{b}\psi\left(
u_{k}^{\prime}\right)  \left|  u_{k}^{\prime\prime}\right|  \,dx\Big)   &
\leq\int_{x_{m}}^{b}\left|  \left(  u^{\prime}\right)  ^{a}\right|
\,dx+\left|  v^{\prime}\right|  \left(  \left]  x_{m},b\right[  \right)\,.
\label{lim3}%
\end{align}
Assume that $\left[  u\right]  \left(  x_{j}\right)
=u_{+}\left(  x_{j}\right)  -u_{-}\left(  x_{j}\right)  >0$. Then \eqref{case1} holds for $k$
sufficiently large.
By (\ref{initial condition II}), \eqref{1inflection}, \eqref{minmax},  and a  change of variables we obtain
\begin{gather}
\int_{x_{j}-\delta_{k}}^{x_{j}+\delta_{k}}\left|  u_{k}^{\prime}\right|
\,dx+\int_{x_{j}-\delta_{k}}^{x_{j}+\delta_{k}}\psi\left(  u_{k}^{\prime
}\right)  \left|  u_{k}^{\prime\prime}\right|  \,dx=\int_{x_{j}-\delta_{k}%
}^{x_{j}+\delta_{k}}|f_{kj}^{\prime}|\,dx+\int_{x_{j}-\delta_{k}}^{x_{j}%
+\delta_{k}}\psi\left(  f_{kj}^{\prime}\right)  \left|  f_{kj}^{\prime\prime
}\right|  \,dx\nonumber\\
\leq f_{kj}\left(  x_{j}+\delta_{k}\right)  -f_{kj}\left(  x_{j}-\delta
_{k}\right)  +\int_{w_{k}\left(  x_{j}-\delta_{k}\right)  }^{f_{kj}^{\prime
}\left(  y_{kj}\right)  }\psi\left(  t\right)  \,dt+\int_{w_{k}\left(
x_{j}+\delta_{k}\right)  }^{f_{kj}^{\prime}\left(  y_{kj}\right)  }\psi\left(
t\right)  \,dt+\frac{2}{k}\,. \label{fkj}%
\end{gather}
By \eqref{1inflection} we have
\begin{equation}
f_{kj}^{\prime}\left(  y_{kj}\right)  =\max_{\left[  x_{j}-\delta_{k}%
,x_{j}+\delta_{k}\right]  }f_{kj}^{\prime}\geq\frac{1}{2\delta_{k}}\left[
f_{kj}\left(  x_{j}+\delta_{k}\right)  -f_{kj}\left(  x_{j}-\delta_{k}\right)
\right]\,.  \label{max fkj}%
\end{equation}
By (\ref{initial condition
I}) and the fact that $w_{k}\rightarrow\left(  u^{\prime}\right)  ^{a}$ in
$L^{1}\left(  \left]  a,b\right[  \right)  $ we obtain
\[
f_{kj}\left(  x_{j}+\delta_{k}\right)  -f_{kj}\left(  x_{j}-\delta_{k}\right)
\rightarrow u_{+}\left(  x_{j}\right)  -\Big(  u_{+}\left(  x_{j-1}\right)
+\int_{x_{j-1}}^{x_{j}}\left(  u^{\prime}\right)  ^{a}\,ds\Big)  =\left[
u\right]  \left(  x_{j}\right)\,.
\]
In turn, using (\ref{max fkj}), we get that $f_{kj}^{\prime}\left(
y_{kj}\right)  \rightarrow\infty$. Thus, letting $k\rightarrow\infty$ in
(\ref{fkj}) and using (\ref{delta k}), we infer
\begin{gather}
\limsup_{k\rightarrow\infty}\Big(  \int_{x_{j}-\delta_{k}}^{x_{j}+\delta_{k}
}\left|  u_{k}^{\prime}\right|  \,dx+\int_{x_{j}-\delta_{k}}^{x_{j}+\delta
_{k}}\psi\left(  u_{k}^{\prime}\right)  \left|  u_{k}^{\prime\prime}\right|
\,dx\Big) \label{lim5}\\
\leq\left[  u\right]  \left(  x_{j}\right)  +\int_{\left(  u^{\prime}\right)
_{-}^{a}\left(  x_{j}\right)  }^{+\infty}\psi\left(  t\right)  \,dt+\int
_{\left(  u^{\prime}\right)  _{+}^{a}\left(  x_{j}\right)  }^{+\infty}%
\psi\left(  t\right)  \,dt\,.\nonumber
\end{gather}
Similarly, if $\left[  u\right]  \left(  x_{j}\right)  =u_{+}\left(
x_{j}\right)  -u_{-}\left(  x_{j}\right)  <0$,  we find
\begin{gather}
\limsup_{k\rightarrow\infty}\Big(  \int_{x_{j}-\delta_{k}}^{x_{j}+\delta_{k}%
}\left|  u_{k}^{\prime}\right|  \,dx+\int_{x_{j}-\delta_{k}}^{x_{j}+\delta
_{k}}\psi\left(  u_{k}^{\prime}\right)  \left|  u_{k}^{\prime\prime}\right|
\,dx\Big) \label{lim4}\\
\leq\left|  \left[  u\right]  \left(  x_{j}\right)  \right|
+\int_{-\infty }^{\left(  u^{\prime}\right)  _{-}^{a}\left(
x_{j}\right)  }\psi\left( t\right)  \,dt+\int_{-\infty}^{\left(
u^{\prime}\right)  _{+}^{a}\left( x_{j}\right)  }\psi\left( t\right)
\,dt\,.\nonumber
\end{gather}
Summing over $j$ in (\ref{lim1}), (\ref{lim5}), (\ref{lim4}) and combining
with (\ref{lim2}), (\ref{lim3}), inequality (\ref{limsup inequality}) follows.

\noindent\textbf{Step 2:} Assume only that $u\in
X^{1}_{\psi}({]a,b[})$ and that $S_{u}$ is finite. We claim that
there exists a sequence $\{u_{k}\}$ such that $u_{k}\rightarrow u$
in $L^{1}\left(  \left]  a,b\right[  \right)$, each $u_{k}$
satisfies the
hypotheses of Step 1, and%
\begin{equation}
\mathcal{G}\left(  u\right)  \geq\limsup_{k\rightarrow\infty}\mathcal{G}%
\left(  u_{k}\right)  . \label{limsup 3}%
\end{equation}

Note that if (\ref{limsup 3}) holds then, by applying Step 1 to each $u_{k}$
we may find a sequence $u_{km}\in W^{2,1}\left(  \left]  a,b\right[  \right)
$ converging to $u_{k}$ in $L^{1}\left(  \left]  a,b\right[  \right)  $\ and
satisfying%
\[
\mathcal{G}\left(  u_{k}\right)  \geq\limsup_{m\rightarrow\infty}%
\mathcal{F}_1\left(  u_{km}\right)  .
\]
By (\ref{limsup 3}) we then have%
\[
\mathcal{G}\left(  u\right)  \geq\limsup_{k\rightarrow\infty}\limsup
_{m\rightarrow\infty}\mathcal{F}_1\left(  u_{km}\right)
\]
and a standard diagonalization argument now yields the existence of
a sequence $m_{k}\to\infty$ such that $u_{km_{k}}\to u$ in
$L^{1}\left(  \left]  a,b\right[  \right)  $\ and
\[
\mathcal{G}\left(  u\right)  \geq\limsup_{k\rightarrow\infty}\mathcal{F}_1
\left(  u_{km_{k}}\right)  .
\]

In the construction of the sequence satisfying (\ref{limsup 3}) we
need to consider the precise representatives
$(u^{\prime})_{\vee}^{a}$ and  $(u^{\prime})_{\wedge}^{a}$ defined
in Remark~\ref{remark limits u'}.  We recall  that
$(u^{\prime})_{\vee }^{a}$ is upper semicontinuous while
$(u^{\prime})_{\wedge}^{a}$ is lower semicontinuous, and so for each
$k\in\mathbb{N}$ we may decompose the open sets $\left\{
(u^{\prime})_{\wedge}^{a}>k\right\}  $ and $\left\{  (u^{\prime
})_{\vee}^{a}<-k\right\}  $ into the union of two finite sequences
of pairwise disjoint open sets $U_{kj}^{+}$ and $U_{kj}^{-}$, that
is,
\[
\bigcup_{j}U_{kj}^{+}=\left\{  (u^{\prime})_{\wedge}^{a}>k\right\}
,\quad\bigcup_{j}U_{kj}^{-}=\left\{  (u^{\prime})_{\vee}^{a}<-k\right\}  \,,
\]
such that
\begin{equation}\label{diam}
{\rm diam\,}(U_{kj}^{+})\leq \Lone(\{  (u^{\prime})_{\wedge}^{a}>k\})\,,\
{\rm diam\,}(U_{kj}^{-})\leq \Lone(\{  (u^{\prime})_{\vee}^{a}<-k\})\quad\text{for every }j\,.
\end{equation}
Note that,
setting  $v_{\vee}:=\Psi_1\circ ( u^{\prime})^{a}_{\vee}$ and
$v_{\wedge}:=\Psi_1\circ ( u^{\prime})^{a}_{\wedge}$, we have
\begin{equation}\label{TV}
|v' |( {]c, d[}  )={\rm Var\,}(  v_{\vee}; {]c,d[ } )={\rm Var\,}(  v_{\wedge};{ ]c,d[ } )
\end{equation}
for every interval ${]c,d[}\subset {]a,b[}$.

For every set $U_{kj}^{\pm}$ we fix a nonnegative function
 $g_{kj}^{\pm}\in C_{c}^{1}(  U_{kj}^{\pm})  $
such that%
\begin{equation}\label{TV33}
\int_{U_{kj}^{\pm}}g_{kj}^{\pm}\left(  x\right)  \,dx=\left(  (u^{\prime}%
)^{c}\right)  ^{\pm}(  U_{kj}^{\pm})  \,,
\end{equation}
and $(  g_{kj}^{\pm})  ^{\prime}$ has
only one zero in the interior of the support of $g_{kj}^{\pm}$. Then we define
\begin{equation}
\label{diversi}
g_{k}^{+}:=\sum\nolimits_j g_{kj}^{+}\,,\quad g_{k}^{-}:=\sum\nolimits_j g_{kj}^{-}\,,
\quad
g_{k}:= g_{k}^{+} - g_{k}^{-}\,,
\quad
w_{k}:=T^{k}_{-k}\circ \left(  u^{\prime}\right)^{a}+g_{k}\,,
\end{equation}
where for any pair of constants $h<k$ the truncation function $T^{k}_{h}$ is defined by
$$
T^{k}_{h}(t):=\begin{cases}
h & \text{for $t\le h\,,$}\\
t & \text{for $h\le t\le k\,,$}\\
k & \text{for $t\ge k\,.$}
\end{cases}
$$
We claim that
\begin{equation}
w_{k}\,\Lone\rightharpoonup\left(  u^{\prime}\right)
^{a}\Lone+\left(  u^{\prime}\right) ^{c}\text{ weakly}^{\ast}\text{
in }M_{b}\left(  \left]  a,b\right[  \right)\,
. \label{weak star conv}%
\end{equation}
Define%
\[
A_{k}:=\left\{  \left(  u^{\prime}\right)  _{\wedge}^{a}>k\right\}
\cup\left\{  \left(  u^{\prime}\right)  _{\vee}^{a}<-k\right\}\,  .
\]
Since by the Chebychev inequality
\begin{equation}
k{\mathcal{L}}^{1}(A_{k})\rightarrow0\,,
\label{Chebicev}%
\end{equation}
it suffices to show that
\begin{equation}
\big(\sum_{j}g_{kj}^{\pm}\big)\,\Lone\rightharpoonup\left(  \left(  u^{\prime}\right)
^{c}\right)  ^{\pm}\text{ weakly}^{\ast}\text{ in }M_{b}\left(  \left]
a,b\right[  \right)  \,. \label{weak star}%
\end{equation}
Let $\varphi\in C_0\left(  \left]  a,b\right[  \right)  $ and $\varepsilon
>0.$ By uniform continuity there exists $\delta=\delta\left(  \varepsilon
\right)  >0$ such that $\left|  \varphi\left(  x\right)  -\varphi\left(
y\right)  \right|  \leq\varepsilon$ for all $x$, $y\in\left]  a,b\right[  $
with $\left|  x-y\right|  \leq\delta$. In view of  \eqref{diam} and (\ref{Chebicev}), for all
$k$ sufficiently large and for all $j$\ we have that ${\rm diam\,}(
U_{kj}^{\pm})  \leq\delta$. Let us fix $y_{kj}^{\pm}\in U_{kj}^{\pm}$. Then,
by \eqref{TV33},%
\begin{align*}
\Big|  \int_{U_{kj}^{\pm}}\varphi\left(  x\right)
g_{kj}^{\pm}&\left( x\right)  \,dx-\int_{U_{kj}^{\pm}}\varphi\left(
x\right)  \,d\left( \left(
u^{\prime}\right)  ^{c}\right)  ^{\pm}\left(  x\right)  \Big|  \\
&=\Big|  \int_{U_{kj}^{\pm}}\left[  \varphi\left(  x\right)
-\varphi( y_{kj}^{\pm})  \right]  g_{kj}^{\pm}\left(  x\right)
\,dx-\int _{U_{kj}^{\pm}}\left[  \varphi\left(  x\right)  -\varphi(
y_{kj}^{\pm })  \right]  \,d\left(  \left(  u^{\prime}\right)
^{c}\right)  ^{\pm
}\left(  x\right)  \Big| \\
&\leq\varepsilon\Big(  \int_{U_{kj}^{\pm}}g_{kj}^{\pm}\left(
x\right) \,dx+\left(  \left(  u^{\prime}\right)  ^{c}\right)
^{\pm}( U_{kj}^{\pm} )  \Big)  \leq2\varepsilon\left(  \left(
u^{\prime}\right)  ^{c}\right)  ^{\pm}(  U_{kj}^{\pm})  \,.
\end{align*}
Summing over $j$ and using the fact that the measures $\big(\sum_{j}g_{kj}^{+}\big)\,\Lone$
and  $\left(  \left(  u^{\prime}\right)  ^{c}\right)  ^{+}$ are concentrated on
$\left\{  \left(  u^{\prime}\right)  _{\wedge}^{a}>k\right\}  $, while the measures
 $\big(\sum_{j}g_{kj}^{-}\big)\,\Lone$ and  $\left(  \left(  u^{\prime}\right)
^{c}\right)  ^{-}$ are concentrated on $\left\{  \left(  u^{\prime}\right)
_{\vee}^{a}<-k\right\}$ (see Remark~\ref{remark limits u'}), we obtain
(\ref{weak star}).

Moreover, we claim that%
\begin{equation}
\lim_{k\rightarrow\infty}\int_{a}^{b}\left|  w_{k}\right|  \,dx=\int_{a}%
^{b}|  \left(  u^{\prime}\right)  ^{a}|  \,dx+|  \left(
u^{\prime}\right)  ^{c}|  \left(  \left]  a,b\right[  \right)\,  .
\label{total var limit}%
\end{equation}
Indeed, using \eqref{TV33}, (\ref{diversi}), and Remark~\ref{remark
limits u'}, we deduce that
\begin{align*}
\int_{a}^{b}\left|  w_{k}\right|  \,dx\leq&\int_{\left\{ |  \left(
u^{\prime}\right)  ^{a}|  \leq k\right\} %
}\left|  \left(  u^{\prime}\right)  ^{a}\right|  \,dx+k{\mathcal{L}}%
^{1}(A_{k})+\sum_{j}\int_{U_{kj}^{+}}g_{kj}^{+}\,dx+\sum_{j}\int_{U_{kj}^{-}}g_{kj}%
^{-}\,dx\\
\leq&\int_{a}^{b}|  \left(  u^{\prime}\right)  ^{a}|  \,dx+k{\mathcal{L}}%
^{1}(A_{k})+\sum
_{j}\left(  \left(  u^{\prime}\right)  ^{c}\right)  ^{+}(  U_{kj}%
^{+})  +\sum_{j}\left(  \left(  u^{\prime}\right)  ^{c}\right)
^{-}(  U_{kj}^{-}) \\
\leq&\int_{a}^{b}|  \left(  u^{\prime}\right)  ^{a}|  \,dx+k{\mathcal{L}}%
^{1}(A_{k})+ |\left(  u^{\prime}\right)  ^{c}|  \left(  \left]  a,b\right[  \right)
\,,
\end{align*}
and the limit superior inequality follows from \eqref{Chebicev}.
The limit inferior inequality follows from (\ref{weak star conv}) and the
lower semicontinuity of the total variation.

Set%
\begin{equation}
u_{k}\left(  x\right)  :=u_{+}\left(  a\right)  +\int_{a}^{x}w_{k}\left(
s\right)  \,ds+\sum_{x_{j}<x,\,x_{j}\in S_{u}}\left[  u\right]  \left(
x_{j}\right)  \label{def uk}%
\end{equation}
and $v_{k}:=\Psi_1\circ (u_{k}^{\prime})^{a}=\Psi_1\circ w_{k}$.

We claim that $u_{k}\rightarrow u$ in $L^{1}\left(  \left]  a,b\right[
\right)  .$ For $x\in\left]  a,b\right[  $ by (\ref{weak star conv}) and
(\ref{total var limit}) it follows that%
\[
\int_{a}^{x}w_{k}\,dy\rightarrow\int_{a}^{x}\left(  u^{\prime}\right)
^{a}\,dy+\left(  u^{\prime}\right)  ^{c}\left(  \left]  a,x\right[  \right)\,,
\]
and so $u_{k}$ converges to $u$ pointwise $\Lone$-a.e.\  and, in turn, in $L^{1}\left(  \left]
a,b\right[  \right)$.

Next we show that
\begin{equation}
\limsup_{k\rightarrow\infty}\mathcal{G}\left(  u_{k}\right)  \leq
\mathcal{G}\left(  u\right)  \,. \label{limsup case 2}%
\end{equation}
From (\ref{total var limit}) we get
\begin{equation}
\left|  u_{k}^{\prime}\right|  \left(  \left]  a,b\right[  \right)
\rightarrow\left|  u^{\prime}\right|  \left(  \left]  a,b\right[
\right) \, .
\label{u'k}%
\end{equation}
Moreover, as $v_k=T^{\Psi_{1}(k)}_{\Psi_{1}(-k)} \circ v$
$\Lone$-a.e.\  in the open set $V_k:={]a,b[}\setminus \supp\,
g_{k}$, we have $\left|  v_{k}^{\prime}\right| \le \left|
v^{\prime}\right| $ as measures in $V_{k}$. In particular, this
yields $\left| v_{k}^{\prime}\right|  \left(  \left]  a,b\right[
\setminus\left( A_{k}\cup S_{u}\right)  \right)  \le \left|
v^{\prime}\right| \left(  \left]  a,b\right[  \setminus\left(
A_{k}\cup S_{u}\right) \right) $ and hence
\begin{equation}
\left|  v_{k}^{\prime}\right|  \left(  \left]  a,b\right[
\setminus\left( A_{k}\cup S_{u}\right)  \right)  \le \left|
v^{\prime}\right|  \left( \left]  a,b\right[  \setminus\left(
A_{\infty}\cup S_{u}\right)  \right) \, ,
\label{v'k less}%
\end{equation}
where
\[
A_{\infty}:=\bigcap_{k}A_{k}=\left\{  \left(  u^{\prime}\right)  _{\wedge}%
^{a}=+\infty\right\}  \cup\left\{  \left(  u^{\prime}\right)  _{\vee}%
^{a}=-\infty\right\}  \,.
\]
Using the properties of $g_{kj}^{+}$ we have%
\begin{gather}
|  v_{k}^{\prime}| (  \left\{  \left(  u^{\prime}\right)
_{\wedge}^{a}>k\right\}  \setminus S_{u})  =\sum_{j}\int_{U_{kj}^{+}%
}\psi(  k+g_{kj}^{+}) | (  g_{kj}^{+})  ^{\prime
}|  \,dx\nonumber\\
=2\sum_{j}\int_{k}^{k+\sup g_{kj}^{+}}\psi(  t)  \,dt\leq
2\mathcal{H}^{0}(
\{  j:\,( (  u^{\prime})  ^{c})  ^{+}(
U_{kj}^{+})  >0\}  )\int_{k}^{\infty}\psi\left(  t\right)  \,dt \,   . \label{v'k}%
\end{gather}
We claim that
\begin{align*}
  2\mathcal{H}^{0}(
\{  j:\,( (  u^{\prime})  ^{c})  ^{+}(
U_{kj}^{+})  >0\}  )&\int_{k}^{\infty}\psi\left(  t\right)  \,dt \\
&  \leq\left|  v^{\prime}\right|  \left(  \left\{  \left(  u^{\prime}\right)
_{\wedge}^{a}>k\right\}  \setminus S_{u}\right)  +4\int_{k}^{\infty}%
\psi\left(  t\right)  \,dt\,.
\end{align*}
Indeed, if   $\left(  \left(  u^{\prime }\right)  ^{c}\right)  ^{+}(
U_{kj}^{+})  >0$, then there exists a connected component
$I_{kj}^{+}=\left]  a_{kj},b_{kj}\right[  $  of $U^{+}_{kj}\setminus
S_{u}$ such that  $\left(  \left(  u^{\prime
}\right)  ^{c}\right)  ^{+}(  I_{kj}^{+})  >0$. Assume that $I_{kj}%
^{+}\subset\subset\left]  a,b\right[  $. Then by Remark~\ref{remark limits u'} we may find
 $c_{kj}\in I_{kj}^{+}$ such that $\left(  u^{\prime
}\right)^{a}_{\wedge}\left(  c_{kj}\right)  =+\infty$, while $\left(  u^{\prime
}\right)  _{\wedge}^{a}\left(  a_{kj}\right)  $, $\left(  u^{\prime}\right)
_{\wedge}^{a}\left(  b_{kj}\right)  \leq k$. Hence by (\ref{TV})
\[
|  v^{\prime}| (  U_{kj}^{+}\setminus S_{u}) \geq |  v^{\prime}| (  I_{kj}^{+})  \geq
2\int_{k}^{\infty}\psi\left(  t\right)  \,dt\,.
\]
Summing over all such intervals and adding the possible contribution of the
intervals $I_{kj}^{+}$ with at least one endpoint in $\left\{  a,b\right\}  $
we obtain the claim. In turn, by (\ref{v'k}) we have
\[
\left|  v_{k}^{\prime}\right|  \left(  \left\{  \left(  u^{\prime}\right)
_{\wedge}^{a}>k\right\}  \setminus S_{u}\right)  \leq\left|  v^{\prime
}\right|  \left(  \left\{  \left(  u^{\prime}\right)  _{\wedge}^{a}>k\right\}
\setminus S_{u}\right)  +4\int_{k}^{\infty}\psi\left(  t\right)  \,dt\,.
\]
A similar estimate holds for the set $\left\{  \left(  u^{\prime}\right)
_{\vee}^{a}<-k\right\}  \setminus S_{u}$ thus yielding%
\begin{equation}
\limsup_{k\rightarrow\infty}\left|  v_{k}^{\prime}\right|  \left(
A_{k}\setminus S_{u}\right)  \leq\left|  v^{\prime}\right|  \left(  A_{\infty
}\setminus S_{u}\right)\,  . \label{v'k bigger}%
\end{equation}
Combining (\ref{v'k less}) with (\ref{v'k bigger}) we obtain%
\begin{equation*}
\limsup_{k\rightarrow\infty}\left|  v_{k}^{\prime}\right|  \left(
\left] a,b\right[  \setminus S_{u}\right)  \leq\left|
v^{\prime}\right|  \left( \left]  a,b\right[  \setminus
S_{u}\right)\,
.%
\end{equation*}
Next we show that%
\begin{equation}
\lim_{k\rightarrow\infty}\sum_{x\in S_{\scriptstyle u_{k}}}\Phi\left(  \nu_{u_{k}},\left(
u_{k}^{\prime}\right)  _{-}^{a},\left(  u_{k}^{\prime}\right)  _{+}%
^{a}\right)  =\sum_{x\in S_{\scriptstyle u}}\Phi\left(  \nu_{u},\left(  u^{\prime}\right)
_{-}^{a},\left(  u^{\prime}\right)  _{+}^{a}\right)  \,. \label{jump}%
\end{equation}
Note that $S_{u_{k}}=S_{u}$ and $\nu_{u_{k}}\left(  x\right)
=\nu_{u}\left( x\right)  $\ for all $k$\ by (\ref{def uk}).\
Moreover, for every $x\in S_{u}$ if $\left(  u^{\prime}\right)
_{+}^{a}\left(  x\right)  \in\mathbb{R}$ then $\left|  \left(
u^{\prime}\right)  _{+}^{a}\left(  y\right)  \right|  \leq k_{0}$
for all $y$ in a right neighborhood of $x$ and for some integer
$k_{0}$.
 Thus,  by (\ref{diversi}) and (\ref{def uk}) we have that $(u_{k}^{\prime}%
)^{a}(y)=\left(  u^{\prime}\right)  ^{a}(y)$ for $k\geq k_{0}$ and for $\Lone$-a.e.\ $y$
in the same right neighborhood. In turn, by \eqref{uprime+} we infer $\left(
u_{k}^{\prime}\right)  _{+}^{a}\left(  x\right)  =\left(  u^{\prime}\right)
_{+}^{a}\left(  x\right)  $ for all $k\geq k_{0}$. If $\left(  u^{\prime
}\right)  _{+}^{a}\left(  x\right)  =\infty$, then for all $k$ we have
 $\left(  u^{\prime}\right)  _{+}^{a} >k$ in a right neighborhood of $x$ by right continuity
(see Remark~\ref{remark limits u'}). By construction this  implies that
$(u_{k}^{\prime})^{a}=w_{k}\geq k$ $\Lone$-a.e.\ in the same right neighborhood. Thus,
$\left(  u_{k}^{\prime}\right)  _{+}^{a}\left(  x\right)  \geq k\rightarrow
\left(  u^{\prime}\right)  _{+}^{a}\left(  x\right)$.
Similarly $\left(  u_{k}^{\prime}\right)  _{-}%
^{a}\left(  x\right)  \rightarrow\left(  u^{\prime}\right)  _{-}^{a}\left(
x\right)  ,$ so that
\[
\Phi\left(  \nu_{u_{k}}\left(  x\right)  ,\left(
u_{k}^{\prime}\right) _{-}^{a}\left(  x\right)  ,\left(
u_{k}^{\prime}\right)  _{+}^{a}\left( x\right)  \right)
\rightarrow\Phi\left(  \nu_{u}\left(  x\right)  ,\left(
u^{\prime}\right)  _{-}^{a}\left(  x\right)  ,\left(
u^{\prime}\right) _{+}^{a}\left(  x\right)  \right) \, .
\]
Hence (\ref{jump}) follows. This, together with (\ref{u'k}) and
(\ref{v'k bigger}), yields (\ref{limsup case 2}).

\noindent\textbf{Step 3:} Let now $u$ be an arbitrary function in $X^1_{\psi
}({]a,b[})$ such that $\mathcal{G}\left(  u\right)  <+\infty$. As
in the previous step it suffices to construct $u_{k}\in X_{\psi}^1({]a,b[})$
satisfying the hypotheses of Step 2, converging to $u$ in $L^{1}\left(
\left]  a,b\right[  \right)  $\ and such that (\ref{limsup case 2}) holds.
Write $S_{u}=\left\{  x_{j}\right\}  $ and for each $k$ define $S_{u}%
^{k}:=\left\{  x_{j}:\,j\leq k\right\}  $ and%
\[
u_{k}\left(  x\right) \: =u_{+}\left(  a\right)  +\int_{a}^{x}\left(
u^{\prime }\right)  ^{a}\,dt+\left(  u^{\prime}\right)  ^{c}\left(
\left]  a,x\right[ \right)  +\sum_{x_{j}<x,\,x_{j}\in
S_{\scriptstyle u}^{\scriptstyle k}}\left[  u\right]  \left(
x_{j}\right)  \,.
\]
It is clear that $\{u_{k}\}$ converges to $u$ in $L^{1}\left( \left]
a,b\right[ \right)  $ and that $\left|  u_{k}^{\prime}\right| \left(
\left] a,b\right[  \right)  \rightarrow\left| u^{\prime}\right|
\left(  \left] a,b\right[  \right)$. Moreover, $\left|
v_{k}^{\prime}\right|  \left( \left]  a,b\right[  \setminus
S_{u}\right)  =\left|  v^{\prime}\right| \left(  \left]  a,b\right[
\setminus S_{u}\right)  $ and
\begin{align*}
\lim_{k\rightarrow\infty}\sum_{x\in S_{\scriptstyle u_{k}}}\Phi\left(  \nu_{u_{k}},\left(
u_{k}^{\prime}\right)  _{-}^{a},\left(  u_{k}^{\prime}\right)  _{+}%
^{a}\right)  &=\lim_{k\rightarrow\infty}\sum_{x\in
S_{u}^{k}}\Phi\left( \nu_{u},\left(  u^{\prime}\right)
_{-}^{a},\left(  u^{\prime}\right)
_{+}^{a}\right)  \\
&=\sum_{x\in S_{\scriptstyle u}}\Phi\left(  \nu_{u},\left(  u^{\prime}\right)  _{-}%
^{a},\left(  u^{\prime}\right)  _{+}^{a}\right) \, .
\end{align*}
This concludes the proof of the theorem.
\end{proof}
We end the section with a compactness result for energy bounded
sequences in $X^1_{\psi}({]a,b[})$.
\begin{corollary}
\label{corollary compactness} Let $\{u_{k}\}$ be a sequence of
functions in $X^1_{\psi}({]a,b[})$ bounded in $L^{1}({]a,b[})$ and
such that
\begin{equation}
C:=\sup_{k}\overline{\mathcal{F}}_1(u_{k})<+\infty\,.
\label{C bound relax energyp=1}%
\end{equation}
Then there exist a subsequence (not relabeled) $\{u_{k}\}$ and a
function $u\in X^1_{\psi}({]a,b[})$ such that
\begin{align}
&  u_{k}\rightharpoonup u\quad\text{weakly}^{\ast}\text{ in }BV({]a,b[}%
)\,,\label{uk weakly up=1}\\
&  \Psi_1\circ(u_{k}^{\prime})^{a}\rightharpoonup\Psi_1\circ(u^{\prime
})^{a}\quad\text{weakly}^{\ast}\text{ in }BV({]a,b[}%
)\,,\label{v'k weakly v'p=1}\\
&
(u_{k}^{\prime})^{a}\rightarrow(u^{\prime})^{a}\quad\text{pointwise
$\Lone$-a.e.\ in }{]a,b[}\,. \nonumber%
\end{align}
\end{corollary}

\begin{proof}
It is well known that convergence in measure is metrizable with the following
metric
\[
d\left(  u_{1},u_{2}\right)  :=\int_{a}^{b}\frac{|u_{1}-u_{2}|}{1+\left|
u_{1}-u_{2}\right|  }\,dx\,,
\]
where $u_{1}$ and $u_{2}$ are (equivalent classes of) measurable functions.

By Theorems~\ref{theorem compactness} and \ref{th:relaxation}, for every $k\in\mathbb{N}$ we may find $w_{k}\in
W^{2,1}({]a,b[})$ such that
\begin{equation}
\int_{a}^{b}|u_{k}-w_{k}|\,dx\leq\frac{1}{k}\,,\quad d\left(  (u_{k}^{\prime
})^{a},w_{k}^{\prime}\right)  \leq\frac{1}{k}\,, \label{diagonal boundsp=1}%
\end{equation}
and
\begin{equation*}
\mathcal{F}_1(w_{k})\leq C+1\,.%
\end{equation*}
By Theorem \ref{theorem compactness} we may find a subsequence (not
relabeled) of $\{w_{k}\}$ and a function $u\in X^1_{\psi}({]a,b[})$
such that (\ref{convw*}), (\ref{wPsi}), (\ref{conv l1}) hold (with
$w_{k}$ in place of $u_{k}$). It now follows from (\ref{diagonal
boundsp=1}) that
$u_{k}\rightarrow u$ in $L^{1}({]a,b[})$ and $(u_{k}^{\prime})^{a}%
\rightarrow(u^{\prime})^{a}$ in measure and hence pointwise
$\mathcal{L}^{1}$ a.e.\ in ${]a,b[}$, up to a further subsequence.
From the bound (\ref{C bound relax energyp=1}), the uniqueness of
the limit, and the invertibility of $\Psi_{1}$, we deduce (\ref{uk
weakly up=1}) and (\ref{v'k weakly v'p=1}).
\end{proof}

\section{The case $p>1$}
\label{s:p>1}

In this section we analyze  the functional \eqref{1dCMM} in the case $p>1$.

Let us state precisely   the standing assumptions.
 Throughout this section $p$ denotes any exponent  in
$]1,+\infty[$,  $\psi\colon{\mathbb{R}}\rightarrow{]0,+\infty\lbrack}$ is a bounded Borel
function satisfying
\begin{equation}
M:=\int_{-\infty}^{+\infty}(\psi(t))^{1/p}\,dt<+\infty\label{intfinp>1}%
\end{equation}
in addition to \eqref{inf>0}, and  $\Psi_{p}\colon\overline
{{\mathbb{R}}}\rightarrow\lbrack0,M]$ denotes the antiderivative of
$\psi^{1/p}$ defined by

\begin{equation}
\Psi_{p}(t):=\int_{-\infty}^{t}(\psi(s))^{1/p}\,ds\,.\label{Psip>1}
\end{equation}
The function $\Psi_{p}^{-1}\colon\lbrack0,M]\rightarrow\overline{{\mathbb{R}}}$
stands for the  inverse function of $\Psi_{p}$.

We now consider the functional ${\mathcal{F}_{p}}\colon L^{1}({]a,b[}%
)\rightarrow\lbrack0,+\infty]$ defined by
\begin{equation}
{\mathcal{F}_{p}}(u):=%
\begin{cases}
\displaystyle  \int_{a}^{b}|u^{\prime}|\,dx+\int_{a}^{b}\psi(u^{\prime
})|u^{\prime\prime}|^{p}\,dx & \text{if }u\in W^{2,p}({]a,b[})\,,\\
+\infty & \text{otherwise.}%
\end{cases}
\label{functional p>1}%
\end{equation}
It turns out that piecewise smooth functions with bounded derivative
and nonempty discontinuity set cannot be approximated by  sequences
with equibounded energy. This is a consequence of Remark
\ref{rm:definition}(i) and Theorem \ref{theorem compactnessp>1}
below, and to this end we introduce a suitable space of functions.
Recall that $Z^{\pm}[(u^{\prime})^{a}]$ are the sets defined in
\eqref{Zu}  and \eqref{Zu-}, while $(u^{\prime})^{s}$ denotes the
singular part of the gradient measure $u^{\prime}$.

\begin{definition}
\label{def spacep>1} Let $X_{\psi}^{p}({]a,b[})$ be the set of all functions
$u\in BV({]a,b[})$ such that $v:=\Psi_{p}\circ(u^{\prime})^{a}$ belongs
to $W^{1,p}({]a,b[})$ and the positive part $\left(  (u^{\prime})^{s}\right)
^{+}$ and the negative part $\left(  (u^{\prime})^{s}\right)  ^{-}$ of the
measure $(u^{\prime})^{s}$ are concentrated on $Z^{+}[(u^{\prime})^{a}]$ and
$Z^{-}[(u^{\prime})^{a}]$, respectively.
\end{definition}

\begin{remark}
\label{rm:definition}

(i) It follows immediately from the definition that if $u\in X_{\psi
}^{p}({]a,b[})$ then $(u^{\prime})^{a}=\Psi_{p}^{-1}(v)$ is
continuous on $[a,b]$ with values in $\overline{\mathbb{R}}$. In
particular, it turns out that
\[
Z^{\pm}[(u^{\prime})^{a}]=\{x\in{]a,b[}:\,(u^{\prime})^{a}=\pm\infty\}\,.
\]
By the assumption on the support of the singular part $(u^{\prime
})^{s}$, we have $\lim_{x\rightarrow x_{0}}(u^{\prime})^{a}(x)=+\infty$ for
every jump point $x_{0}$ with $u_{+}(x_{0})-u_{-}(x_{0})>0$ and
$\lim_{x\rightarrow x_{0}}(u^{\prime})^{a}(x)=-\infty$ for every jump point
$x_{0}$ with $u_{+}(x_{0})-u_{-}(x_{0})<0$. This means that if
$S_{u}$ is nonempty then $u$ cannot have bounded derivative outside the
discontinuity set. In particular, piecewise constant functions are not included in the class
$X_{\psi}^{p}({]a,b[})$.

(ii) We observe that the function $(u^{\prime})^{a}$ is differentiable
$\mathcal{L}^{1}$ a.e.\ in ${]a,b[}$ with%
\begin{equation}
v^{\prime}=\psi^{\frac{1}{p}}\left(  (u^{\prime})^{a}\right)  \left(
(u^{\prime})^{a}\right)  ^{\prime}\,. \label{v' equal blah}%
\end{equation}
To see this, we consider the open set
\[
A_{k}:=\{x\in{]a,b[}:\,-k<(u^{\prime})^{a}<k\}\,.
\]
Since by (\ref{inf>0}) the function $\Psi_{p}^{-1}$ is Lipschitz
continuous in the interval $[\Psi_{p}(-k), \Psi_{p}(k)]$ and $v\in
W^{1,p}({]a,b[})$, by the chain rule we have that
$(u^{\prime})^{a}=\Psi_{p}^{-1}\circ v \in W^{1,p}(A_{k})$ and, in
particular, it is differentiable $\mathcal{L}^{1}$-a.e.\ in $A_{k}$
and (\ref{v' equal blah})
holds. Since $(u^{\prime})^{a}$ is integrable we have that%
\[
\mathcal{L}^{1}\Big(  {]a,b[}\setminus\bigcup_{k}A_{k}\Big)  =0
\]
and the conclusion follows.

(iii) It is easy to check that $X_{\psi}^{p}({]a,b[})$ may contain
discontinuous functions. An example is given by the following construction:
Let $\psi\colon{\mathbb{R}}\rightarrow{]0,+\infty\lbrack}$ be defined by
\[
\psi(t):=%
\begin{cases}
1 & \text{if }|t|\leq1\,,\\
\frac{1}{|t|^{\alpha}} & \text{if }|t|>1\,,
\end{cases}
\]
where $\alpha$ is any number in ${]1,+\infty[}$, and let $p\in {]1,\frac{\alpha+1}2[}$.
Consider now the discontinuous
functions $u\colon {]{-1},1[}\rightarrow{\mathbb{R}}$ given by
\[
u(x):=%
\begin{cases}
-|x|^{\beta} & \text{if }x\leq0\,,\\
1+x^{\beta} & \text{if }x>0\,,
\end{cases}
\]
with
$$
0<\beta<1-\frac{p-1}{\alpha-p}\,.
$$
A straightforward computation shows that the function $\Psi_{p}\circ(u^{\prime})^{a}$ belongs to $W^{1,p}%
(]{-}1,1[)$, which in turn implies that $u\in X_{\psi}^{p}(]{-}1,1[)$.

(iv)  Finally,  the same construction of Proposition~\ref{prop:cantor} shows that for
every admissible $\psi$ satisfying \eqref{algebraic}  the space $X_{\psi}%
^{p}({]a,b[})$ contains a function with nontrivial Cantor part,  if  $p$
is sufficiently close to $1$. We omit  the details of this fact which can be easily checked following
 step by step the proof of  Proposition~\ref{prop:cantor}.
\end{remark}

The next theorem is the counterpart of Theorem~\ref{theorem
compactness} for the case $p>1$. It establishes that energy bounded
sequences  are relatively compact in $X_{\psi}^{p}({]a,b[})$. The
proof is similar to the one of Theorem~\ref{theorem compactness},
nevertheless since this is the main result of this section we
reproduce it in full detail for the reader's convenience.

\begin{theorem}
\label{theorem compactnessp>1} Let $\{u_{k}\}$ be a sequence of
functions bounded in $L^{1}({]a,b[})$ and such that
\begin{equation}
C:=\sup_{k}{\mathcal{F}_{p}}(u_{k})<+\infty\,.\label{supenergyp>1}%
\end{equation}
Then there exist a subsequence (not relabeled) $\{u_{k}\}$ and a
function $u\in X_{\psi}^{p}({]a,b[})$ such that
\begin{align}
&  u_{k}\rightharpoonup u\quad\text{weakly}^{\ast}\text{ in }BV({]a,b[}%
)\,,\label{convw*p>1}\\
&  \Psi_{p}\circ u_{k}^{\prime}\rightharpoonup\Psi_{p}\circ(u^{\prime}%
)^{a}\quad\text{weakly in }W^{1,p}({]a,b[})\,,\nonumber\\
&  u_{k}^{\prime}\rightarrow(u^{\prime})^{a}\quad\text{pointwise in
}{]a,b[}\,.\label{conv l1p>1}%
\end{align}
\end{theorem}

\begin{proof}
By (\ref{functional p>1}) and (\ref{supenergyp>1}) we may assume that each $u_{k}$
belongs to in $W^{2,p}({]a,b[})$ and that
\begin{equation}
C_{1}:=\sup_{k}\int_{a}^{b}[\,|u_{k}|+|u_{k}^{\prime}|+\psi(u_{k}^{\prime
})|u_{k}^{\prime\prime}|^{p}\,]\,dx<+\infty\,. \label{BoundCp>1}%
\end{equation}
Let us define
\begin{equation}
v_{k}:=\Psi_{p}\circ u_{k}^{\prime}\,. \label{vkp>1}%
\end{equation}
As $\Psi_{p}$ is Lipschitz in ${\mathbb{R}}$, the functions $v_{k}$
belong to $W^{1,p}({]a,b[})$ and
\begin{equation}
v_{k}^{\prime}=(\psi(u_{k}^{\prime}))^{1/p}u_{k}^{\prime\prime}%
\qquad{\mathcal{L}}^{1}\text{-a.e.\ on }{]a,b[}\,. \label{derivvkp>1}%
\end{equation}
It follows from (\ref{intfinp>1}) and (\ref{supenergyp>1}) that
\begin{equation}
\int_{a}^{b}[\,|v_{k}|^{p}+|v_{k}^{\prime}|^{p}\,]\,dx\leq M^{p}(b-a)+C_{1}\,.
\label{BoundC2p>1}%
\end{equation}
By (\ref{BoundCp>1}) and (\ref{BoundC2p>1}), passing to a subsequence (not
relabeled), we may assume that
\[
u_{k}\rightharpoonup u\quad\text{weakly}^{\ast}\text{ in }BV({]a,b[})
\]
and
\begin{equation}
v_{k}\rightharpoonup v \quad\text{weakly in } W^{1,p}({]a,b[})
\label{convpoint3p>1}%
\end{equation}
for some functions $u\in BV({]a,b[})$ and $v\in W^{1,p}({]a,b[}; [0, M])$.

Since $\Psi_{p}^{-1}$ is continuous, we obtain
\begin{equation}
u_{k}^{\prime}=\Psi_{p}^{-1}\circ v_{k}\rightarrow w:=\Psi_{p}^{-1}\circ
v\quad\text{pointwise in }{]a,b[}\,. \label{convpoint4p>1}%
\end{equation}
Note also that $w$ is continuous with values in $\overline{\mathbb{R}}$.

\noindent We now split the remaining part of the proof into two steps.

\noindent\textbf{Step 1: }We prove that
\begin{equation}
w=(u^{\prime})^{a}\quad\Lone\text{-a.e.\ on }{]a,b[}\,. \label{w=dotup>1}%
\end{equation}
If not, arguing as for \eqref{wneqdotu}, we may find $t_{0}>0$ and an
infinite number of disjoint open intervals $I$ such that
\begin{equation}
{\mathcal{L}}^{1}(\{w\neq(u^{\prime})^{a}\}\cap\{|w|<t_{0}\}\cap I)>0\,.
\label{wneqdotup>1}%
\end{equation}
By H\"older's inequality and a change of variables we obtain
\begin{equation}
\begin{array}{c}
\displaystyle \int_{I}\psi(u_{k}^{\prime})|u_{k}^{\prime\prime}|^{p}\,dx\geq\frac{1}%
{{\mathcal{L}}^{1}(I)^{p-1}} \Big(  \int_{I}(\psi(u_{k}^{\prime}%
))^{1/p}|u_{k}^{\prime\prime}|\,dx\Big)  ^{p}\smallskip\\
\displaystyle\,\quad \qquad\qquad\qquad\geq \frac{1}{(b-a)^{p-1}
}\Big( \int_{m_{k}}^{M_{k}}(\psi(t))^{1/p}\,dt\Big)  ^{p}\,,
\end{array}
\label{changevar2p>1}
\end{equation}
where $m_{k}:=\inf_{I}u_{k}^{\prime}$ and
$M_{k}:=\sup_{I}u_{k}^{\prime}$.

Reasoning as in the first step of the proof of Theorem~\ref{theorem
compactness}, we can show that at least one of the two sequences
$\{m_{k}\}$ and $\{M_{k}\}$ is divergent. If
$\lim_{k}\,M_{k}=+\infty$
then by \eqref{convpoint4p>1}  $\limsup_{k}%
\,m_{k}<t_{0}$ and, in turn, from \eqref{changevar2p>1}  we obtain
\[
\liminf_{k\rightarrow\infty}\int_{I}\psi(u_{k}^{\prime})|u_{k}^{\prime\prime
}|^{p}\,dx\geq\frac{1}{(b-a)^{p-1}}\Big(  \int_{t_{0}}^{+\infty}%
(\psi(t))^{1/p}\,dt\Big)  ^{p}>0\,.
\]
Analogously, if $\lim_{k}\,m_{k}=-\infty$ then
\[
\liminf_{k\rightarrow\infty}\int_{I}\psi(u_{k}^{\prime})|u_{k}^{\prime\prime
}|^{p}\,dx\geq\frac{1}{(b-a)^{p-1}}\Big(  \int_{-\infty}^{t_{0}}%
(\psi(t))^{1/p}\,dt\Big)  ^{p}>0\,.
\]
In any case for an arbitrarily large number $m$ of disjoint intervals $I$
satisfying \eqref{wneqdotup>1}  , adding the contributions of each interval we
obtain
\[
\liminf_{k\rightarrow\infty}\int_{a}^{b}\psi(u_{k}^{\prime})|u_{k}%
^{\prime\prime}|^{p}\,dx\geq\frac{m}{(b-a)^{p-1}}\min\Big\{  \Big(
\int_{t_{0}}^{+\infty}(\psi(t))^{1/p}\,dt\Big)  ^{p}, \Big(  \int_{-\infty
}^{t_{0}}(\psi(t))^{1/p}\,dt\Big)  ^{p}\Big\}  \,,
\]
which contradicts~\eqref{BoundCp>1}  for $m$ large enough. This concludes the
proof of~\eqref{w=dotup>1}  and, in turn, of \eqref{conv l1p>1}  .

\noindent\textbf{Step 2:} To prove that $u\in X_{\psi}^{p}({]a,b[})$ it remains
to show that the positive part $\left(  (u^{\prime})^{s}\right)  ^{+}$ and the
negative part $\left(  (u^{\prime})^{s}\right)  ^{-}$ of the measure
$(u^{\prime})^{s}$ are concentrated on $Z^{+}[(u^{\prime})^{a}]$ and
$Z^{-}[(u^{\prime})^{a}]$ respectively.

Arguing as in Step 2 of the proof of Theorem~\ref{theorem
compactness}, one can see that it is enough to show
\begin{equation}
E^{+}[u^{\prime}]\setminus Z^{+}[(u^{\prime})^{a}]\text{ and }E^{-}[u^{\prime
}]\setminus Z^{-}[(u^{\prime})^{a}]\text{ are empty,} \label{countablep>1}%
\end{equation}
where $E^{+}[u^{\prime}]$ and $E^{-}[u^{\prime}]$ are the sets introduced in
\eqref{Eu+}  and \eqref{Eu-}  . We only show that $E^{+}[u^{\prime}]\setminus
Z^{+}[(u^{\prime})^{a}]$ is empty, since the other property can be proved in
the same way.

Assume by contradiction that $E^{+}[u^{\prime}]\setminus Z^{+}[(u^{\prime
})^{a}]$ contains a point $x_{0}$. Denote $t_{0}:=2|(w(x_{0})|$, fix any
 $t_{1}>t_{0}$, and choose $\varepsilon_{0}>0$ such that
\begin{equation}
\frac{1}{(2\varepsilon_{0})^{p-1}}\Big(  \int_{t_{0}}^{t_{1}}(\psi
(t))^{1/p}\,dt\Big)  ^{p}>C\,, \label{mintp>1}%
\end{equation}
where $C$ is the constant appearing in~(\ref{supenergyp>1}). By \eqref{Eu+}
there exists $0<\varepsilon<\varepsilon_{0}$ such that
\begin{equation}
\frac{\left(  u^{\prime}\right)  ^{+}({]x_{0}-\varepsilon,x_{0}+\varepsilon
\lbrack})}{2\varepsilon}>t_{1}\,. \label{quotientp>1}%
\end{equation}
Set $I:={]x_{0}-\varepsilon,x_{0}+\varepsilon\lbrack}$. By
H\"older's inequality and a change of variables  (see
\eqref{changevar2p>1}) we obtain
\begin{equation}
\int_{I}\psi(u_{k}^{\prime})|u_{k}^{\prime\prime}|^{p}\,dx\geq\frac{1}%
{(2\varepsilon_{0})^{p-1}}\Big(  \int_{m_{k}}^{M_{k}}(\psi(t))^{1/p}%
\,dt\Big)  ^{p}\,, \label{changevar10p>1}%
\end{equation}
where $m_{k}:=\inf_{I}u_{k}^{\prime}$ and
$M_{k}:=\sup_{I}u_{k}^{\prime}$. By \eqref{convpoint4p>1}  and the
fact that $w(x_{0}) <t_{0}$, we deduce that
\begin{equation}
\limsup_{k\rightarrow\infty}\,m_{k}<t_{0}\,. \label{inf20p>1}%
\end{equation}
On the other hand, reasoning as at the end of the proof of Theorem~\ref{theorem compactness}, we deduce from \eqref{convw*p>1} and \eqref{quotientp>1} that
$$
\liminf_{k\rightarrow\infty}\frac{1}{2\varepsilon}\int_{x_{0}-\varepsilon
}^{x_{0}+\varepsilon}(u_{k}^{\prime})^{+}\,dx
   \geq \frac{\left(  u^{\prime}\right)  ^{+}({]x_{0}-\varepsilon,x_{0}%
+\varepsilon\lbrack})}{2\varepsilon}>t_{1}\,,
$$
which implies that
\begin{equation}
\liminf_{k\rightarrow\infty}\,M_{k}>t_{1}\,. \label{sup20p>1}%
\end{equation}
{}From \eqref{mintp>1}, \eqref{changevar10p>1}, \eqref{inf20p>1}, and
\eqref{sup20p>1}  we obtain
\[
\liminf_{k\rightarrow\infty}\int_{I}\psi(u_{k}^{\prime})|u_{k}^{\prime\prime
}|^{p}\,dx\geq\frac{1}{(2\varepsilon_{0})^{p-1}}\Big(  \int_{t_{0}}^{t_{1}%
}(\psi(t))^{1/p}\,dt\Big)  ^{p}>C\,,
\]
which contradicts~\eqref{BoundCp>1}  . This shows (\ref{countablep>1}) and
concludes the proof of the theorem.
\end{proof}

We next identify the relaxation  of $\mathcal{F}_{p}$ with respect
to strong convergence in $L^{1}({]a,b[})$.
\begin{theorem}
Let $\overline{\mathcal{F}}_{p}\colon
L^{1}({]a,b[})\rightarrow\lbrack0,+\infty]$ be defined by %
\begin{equation}
\overline{\mathcal{F}}_{p}\left(  u\right)  :=\inf\left\{
 \liminf_{k\rightarrow\infty}\mathcal{F}_{p}\left(  u_{k}\right)  :\,u_{k}\rightarrow
u\text{ in }L^{1}({]a,b[})\right\}  \label{relaxed functional p>1}%
\end{equation}
for every $u\in L^{1}({]a,b[})$.
 Then
\begin{equation}
\overline{\mathcal{F}}_{p}\left(  u\right)  =%
\begin{cases}
\displaystyle  |u^{\prime}|({]a,b[})+\int_{a}^{b}|v^{\prime}|^{p}\,dx\, &
\text{if }u\in X_{\psi}^{p}({]a,b[})\,,\\
+\infty & \text{otherwise,}%
\end{cases}
\label{relaxed representation p>1}%
\end{equation}
where $v:=\Psi_{p}\circ(u^{\prime})^{a}$.
\end{theorem}

\begin{proof}
We  sketch the proof focusing only on the main changes with respect to the proof of Theorem~\ref{th:relaxation}.
Let $\mathcal{G}_p$ be the functional defined by the right hand side
of (\ref{relaxed representation p>1}).\par We start by showing that
\begin{equation}
\mathcal{G}_p(u)\leq\liminf_{k\rightarrow\infty}{\mathcal{F}}_p(u_{k}%
)\,.\label{liminfp>1}%
\end{equation}
whenever $u_{k}\rightarrow u$ in $L^{1}({]a,b[})$. It is enough to
consider sequences $\{u_{k}\}$ for which the liminf is a limit and
has a finite value. Then $u_{k}$ belongs to $W^{2,1}({]a,b[})$ and
 \eqref{supenergyp>1}  is satisfied. Setting $v_k:=\Psi_p\circ u'_k$, by Theorem~\ref{theorem compactnessp>1} we
 have $v_k\rightharpoonup v$ weakly in $W^{1,p}({]a,b[})$. Using the fact that  $|v'_k|^p=\psi(u'_k)|u_k''|^p$, we
 deduce that
\begin{equation}\label{vkweakly}
 \int_a^b|v'|^p\, dx\leq \liminf_{k\to\infty}\int_a^b\psi(u'_k)|u_k''|^p\, dx\,.
\end{equation}
Inequality \eqref{liminfp>1}  follows now from \eqref{vkweakly} and the lower semicontinuity of the
total variation.

We split the proof of the limsup inequality into several steps.

\noindent {\bf Step 1:} Let $u\in  X_{\psi}^p({]a,b[})$ be such that
$(u')^s=0$. We claim that there exists a sequence $\{u_k\}$ in
$W^{2,p}({]a,b[})$ such that $u_k\to u$ in $L^{1}({]a,b[})$ and
\begin{equation}
\limsup_{k\rightarrow\infty}{\mathcal{F}}_p(u_{k})\leq \mathcal{G}_p(u)\,.
\label{limsupp>1}%
\end{equation}
Define $w_k:=((u')^a\lor -k)\land k$. Using the fact that $(u')^a\in W^{1,p}(A_{2k})$, where
$$
A_{2k}:=\{x\in{]a,b[}:\,-2k<(u^{\prime})^{a}<2k\}
$$
as observed in Remark~\ref{rm:definition}-(ii), one sees that $w_k\in W^{1,p}({]a,b[})$.
Define
$$
u_k(x):=u_+(a)+\int_a^xw_k(y)\, dy\,.
$$
It is easy to see that  $u_k\to u$ in $L^{1}({]a,b[})$ and  \eqref{limsupp>1} holds.

\noindent {\bf Step 2:} Assume that $u\in  X_{\psi}^p({]a,b[})$, $(u')^c=0$, and $S_u$ is finite. We claim that there exists
 a sequence $\{u_k\}$ of functions in $X_{\psi}^p({]a,b[})$, with $(u'_k)^s=0$, such that $u_k\to u$ in $L^{1}({]a,b[})$
and
\begin{equation}
\limsup_{k\rightarrow\infty}{\mathcal{G}}_p(u_{k})\leq \mathcal{G}_p(u)\,.
\label{limsup2p>1}%
\end{equation}
Since the construction is local, it is enough to consider the case
$S_u=\{x_0\}$ for some $x_0\in{]a,b[}$ with $[u](x_0)>0$. By the
properties of $X_\psi^p({]a,b[})$ we can find two sequences
$x_k\nearrow x_0$ and $y_k\searrow x_0$ such that
$$
u(x_k)\to u_-(x_0)\,,\quad u(y_k)\to u_+(x_0)\,,\quad\text{and}\quad
 (u')^a(x_k)=(u')^a(y_k)\to
(u')^a(x_0)=+\infty\,.
$$
Consider the affine functions $h_k(x):=u(x_k)+(u')^a(x_k)(x-x_k)$.
For every $k$ sufficiently large there exists $z_k\in{]x_k,b[}$ such
that $h_k(z_k)=u_+(x_0)$. Since $(u')^a(x_k)\to+\infty$ and $x_k\to
x_0$, we have that $z_k\to x_0$ as $k\to\infty$. Define
$$
u_k(x):=\begin{cases}
u(x) & \text{if $a<x\leq x_k\,,$}\\
h_k(x) & \text{if $x_k<x\leq z_k\,,$}\\
u(x+y_k-z_k)+u_+(x_0)-u(y_k) & \text{if $z_k<x<b\,.$}
\end{cases}
$$
Using the fact that $(u')^a(x_k)=(u')^a(y_k)$, it is easy to check
that $u_k\in X_{\psi}^p({]a,b[})$, with $(u'_k)^s=0$,  $u_k\to u$ in
$L^{1}({]a,b[})$, and \eqref{limsup2p>1} holds.

\noindent {\bf Step 3:} Assume that $u\in  X_{\psi}^p({]a,b[})$ and $(u')^c=0$. We claim that there exists
 a sequence of functions $u_k$ in $X_{\psi}^p({]a,b[})$, with $(u'_k)^c=0$ and $S_{u_k}$  finite, such that $u_k\to u$
 in $L^{1}({]a,b[})$ and \eqref{limsup2p>1} holds.

To see this, it is enough to consider the same approximation constructed in Step 3 of the proof of Theorem~\ref{th:relaxation}.

\noindent {\bf Step 4:} Assume that $u\in  X_{\psi}^p({]a,b[})$. We claim that there exists
 a sequence of functions $u_k$ in $X_{\psi}^p({]a,b[})$, with $(u'_k)^c=0$, such that $u_k\to u$
 in $L^{1}({]a,b[})$ and \eqref{limsup2p>1} holds.

Since $(u')^a$ is continuous from ${]a,b[}$ into $\overline\R$ and integrable
(see Remark~\ref{rm:definition}), we have
 that $K:=\{x\in{]a,b[}:\, |(u')^a|=+\infty\}$ is relatively closed in ${]a,b[}$ with zero
$\Lone$ measure. Hence, we may find a sequence
of open sets $A_k\subset {]a,b[}$  such that $A_k\searrow K$. Let $\{I^k_{j}\}_j$ be the collection of all connected components of $A_k$ intersecting $K$.
Let $c^k_j:=(u')^s(I^k_{j})>0$. By the properties of $X_{\psi}^p({]a,b[})$ for every $j$
we may choose $x^k_j\in I^k_{j}\cap K$ such that $(u')^a(x^k_j)=+\infty$ if $c^k_j>0$ and
$(u')^a(x^k_j)=-\infty$ if $c^k_j<0$.
Define
$$
u_k(x):=u_+(a)+\int_a^x (u')^a(y)\, dy+ \sum_{j:\, x^k_j\leq
x}c^k_j\,.
$$
Using the definition of $X_{\psi}^p({]a,b[})$ one can check that
$$
\sum_j c^k_j\delta_{x^k_j}\rightharpoonup (u')^s\quad\text{weakly}^* \text{ in } M_b({]a,b[})
$$
and $|\sum_j c^k_j\delta_{x^k_j}|({]a,b[})\to |(u')^s|({]a,b[})$ as
$k\to\infty$. Using this fact it is easy to see that the sequence
$\{u_k\}$ meets all the requirements.

By combining Steps 1-4 with a diagonal argument one can finally prove that \eqref{limsupp>1} holds for every $u$ in
$X_{\psi}^p({]a,b[})$.
\end{proof}

\begin{corollary}
\label{corollary compactness p>1} Let $\{u_{k}\}$ be a sequence of
functions in $X_{\psi}^{p}({]a,b[})$ bounded in $L^{1}({]a,b[})$ and
such that
\begin{equation}
C:=\sup_{k}\overline{\mathcal{F}}_{p}(u_{k})<+\infty\,.
\label{C bound relax energy}%
\end{equation}
Then there exists a subsequence (not relabeled) $\{u_{k}\}$ and a
function $u\in X_{\psi}^{p}({]a,b[})$ such that
\begin{align}
&  u_{k}\rightharpoonup u\quad\text{weakly}^{\ast}\text{ in }BV({]a,b[}%
)\,,\label{uk weakly u}\\
&  \Psi_{p}\circ(u_{k}^{\prime})^{a}\rightharpoonup\Psi_{p}\circ(u^{\prime
})^{a}\quad\text{weakly in }W^{1,p}({]a,b[})\,,\label{v'k weakly v'}\\
&  (u_{k}^{\prime})^{a}\rightarrow(u^{\prime})^{a}\quad\text{pointwise in
}{]a,b[}\,. \label{u'ka pointwise u'a}%
\end{align}
\end{corollary}

\begin{proof}

With an argument entirely similar to the one used in the proof of Corollary~\ref{corollary compactness} we can
 extract a subsequence $\{u_k\}$ which satisfies
 (\ref{uk weakly u}) and (\ref{v'k weakly v'}). In turn
(\ref{v'k weakly v'}) and the continuity of $\Psi_{p}^{-1}$\ in $\overline
{\mathbb{R}}$ imply (\ref{u'ka pointwise u'a}).
\end{proof}

\section{The staircase effect}
\label{s:staircase}

The purpose of this section is to show analytically that the presence of the
higher order term in the functional $\overline{\mathcal{F}}$ prevents the
occurrence of the so-called \emph{staircase effect} as opposed to what happens
in image reconstructions based on the total variation functional.

\subsection{The Rudin-Osher-Fatemi model}
We start by showing that staircase-like structures do appear in
solutions to the  Rudin-Osher-Fatemi problem; i.e., in minimizers for the functional
$\mathrm{ROF}_{\lambda,g}:BV({]a,b[})\to \R$ defined by
\[
\mathrm{ROF}_{\lambda,g}(w):=|w^{\prime}|({]a,b[})+\lambda
\int_{a}^{b}(w-g)^{2}\,dx\,,
\]
where $\lambda>0$ is the \emph{fidelity parameter} and $g\in L^{2}({]a,b[})$ is
the given ``signal''\ to be processed. This fact is well known and numerically
observed in many situations. We provide here
a simple analytical example. A different example can be found in~\cite{Ri}.
It will be constructed by means of the
following proposition which deals with minimizers of $\mathrm{ROF}_{\lambda
,g}$ when $g$ is a monotone function.

\begin{proposition}
\label{rof}Let $g:[a,b]\rightarrow\lbrack0,1]$ be a nondecreasing function
such that $g_+(a)=0$ and $g_-(b)=1$. Let $g^{-1}$ denote the left-continuous
generalized inverse
of $g$, defined by
\begin{equation}
g^{-1}(c):=\inf\{x\in\lbrack a,b]:\,g(x)\geq c\} \label{gif}%
\end{equation}
for every $c\in\lbrack0,1]$ and assume that there exist $0<c_{1}<c_{2}<1$ such
that
\begin{equation}
2\lambda\int_{a}^{g^{-1}(c_{1})}(c_{1}-g(x))\,dx=1\qquad\text{and}\qquad
2\lambda\int_{g^{-1}(c_{2})}^{b}(g(x)-c_{2})\,dx=1\,. \label{c1c2}%
\end{equation}
Then the function $u$, defined by
\begin{equation*}
u(x):=%
\begin{cases}
c_{1} & \text{if }a\leq x\leq g^{-1}(c_{1})\,,\\
g(x) & \text{if }g^{-1}(c_{1})<x\leq g^{-1}(c_{2})\,,\\
c_{2} & \text{if }g^{-1}(c_{2})<x\leq b\,,
\end{cases}
\end{equation*}
is the unique minimizer of $\mathrm{ROF}_{\lambda,g}$ in $BV({]a,b[})$.
\end{proposition}

\begin{remark}
\label{rm:c1c2} Since
$$
\int_{a}^{g^{-1}(c)}(c-g(x))\,dx=\int_0^{c}g^{-1}(y)\, dy\,,
\qquad
\int_{g^{-1}(c)}^{b}(g(x)-c)\,dx= \int_{c}^1g^{-1}(y)\, dy
$$
for all $c\in[0,1]$, the continuity of the integral implies that condition \eqref{c1c2}
is satisfied for
every $\lambda$ sufficiently large.
\end{remark}

\begin{proof}
[Proof of Proposition~\ref{rof}]We split the proof into two steps.

\noindent\textbf{Step 1.} We assume first that  $u$
 is absolutely continuous.
 In order to prove the minimality of $u$, by density it suffices to
show that $\mathrm{ROF}_{\lambda,
g}(u+\varphi)\geq\mathrm{ROF}_{\lambda, g}(u)$ for every $\varphi\in
C^{1}([a,b])$, which, in turn, due to the convexity of
$\mathrm{ROF}_{\lambda, g}$, is equivalent to proving that
\begin{equation}
\label{suff1}\left.  \frac{d^+}{d\varepsilon} \mathrm{ROF}_{\lambda,
g}(u+\varepsilon\varphi)\right|  _{\varepsilon=0} \geq0 \quad\text{for every }
\varphi\in C^{1}([a,b])\,,
\end{equation}
where $\frac{d^+}{d\varepsilon}$ denotes the right derivative.
By a straightforward computation we have
\begin{equation}
\label{dde}\left.  \frac{d^+}{d\varepsilon} \mathrm{ROF}_{\lambda,
g}(u+\varepsilon\varphi)\right|  _{\varepsilon=0} =\int_{\{u^{\prime}%
=0\}}|\varphi^{\prime}|\,
dx+\int_{\{u^{\prime}>0\}}\varphi^{\prime}\,
dx+2\lambda\int_{a}^{b}(u-g)\varphi\, dx\,.
\end{equation}
Consider now the function $\theta:[a,b]\to[0,1]$ defined by $\theta
(x):=2\lambda\int_{a}^{x} (u-g)\, dt$. Using \eqref{c1c2}  and the definition
of $u$ one can check that $\theta(a)=\theta(b)=0$, $0\leq\theta\leq1$, and
$\theta\equiv1$ in $[g^{-1}(c_{1}), g^{-1}(c_{2})]$. In particular,
$\{u^{\prime}>0\}\subset[g^{-1}(c_{1}), g^{-1}(c_{2})]\subset\{\theta=1\}$ so
that by \eqref{dde}
\[
\left.  \frac{d^+}{d\varepsilon} \mathrm{ROF}_{\lambda, g}(u+\varepsilon
\varphi)\right|  _{\varepsilon=0} \geq\int_{a}^{b}\varphi^{\prime}\theta\,
dx+2\lambda\int_{a}^{b}(u-g)\varphi\, dx=0\,,
\]
where the last equality is obtained by integrating by parts and by using the
fact that $\theta^{\prime}=2\lambda(u-g)$ and $\theta(a)=\theta(b)=0$. This
shows \eqref{suff1}  and concludes the proof of Step 1.

\noindent\textbf{Step 2.} In the general case, we construct a
sequence $\{g_{k}\}\subset AC([g^{-1}(c_{1}), g^{-1}(c_{2})])$ of
nondecreasing functions such that $g_{k}(g^{-1}(c_{1}))=c_{1}$,
$g_{k}(g^{-1}(c_{2}))=c_{2}$, and $g_{k}\to g$ in
$L^{2}([g^{-1}(c_{1}), g^{-1}(c_{2})])$. Let $\tilde g_{k}$ be the
function that coincides with $g_{k}$ in $[g^{-1}(c_{1}),
g^{-1}(c_{2})]$ and with $g$ elsewhere in $[a,b]$ and, analogously,
set $u_{k}$ to be equal to $g_{k}$ in $[g^{-1}(c_{1}),
g^{-1}(c_{2})]$ and to $u$ elsewhere. For any $v\in BV({]a,b[})$, by
applying the previous step we obtain
\[
\mathrm{ROF}_{\lambda, \tilde g_{k}}(v)\geq\mathrm{ROF}_{\lambda, \tilde
g_{k}}(u_{k})= \mathrm{ROF}_{\lambda, g}(u)\,.
\]
The minimality of $u$ follows by letting $k\to\infty$. Finally, uniqueness is
a consequence of the strict convexity of $\mathrm{ROF}_{\lambda, g}$.
\end{proof}

As a corollary of the previous result we can prove analytically the occurrence
of the staircase effect in a very simple case. Let $g\left(  x\right)  :=x$,
$x\in\left[  0,1\right]  $, be the original $1$D image to which we add the
"noise"%
\[
h_{n}\left(  x\right)  :=\frac{i}{n}-x\text{\qquad if\ \
}\frac{i-1}{n}\leq x<\frac{i}{n}\,\text{,\ \ }i=1,\ldots,n\,,%
\]
where $n\in\mathbb{N}$, so that the resulting degraded $1$D image is given by
the staircase function
\begin{equation}
g_{n}\left(  x\right)  :=\frac{i}{n}\text{\qquad if\ \ }\frac{i-1}{n}\leq
x<\frac{i}{n}\,\text{,\ \ }i=1\text{, \ldots, }n\text{.} \label{starcase functions}%
\end{equation}
Note that, even though $h_{n}\rightarrow0$ uniformly, the reconstructed image
$u_{n}$ preserves the staircase structure of $g_{n}$. Indeed, we show that
there exists a non degenerate interval $I\subset\left[  0,1\right]  $ such
that each $u_{n}$ coincides with the degraded $1$D image $g_{n}$ in $I$ for
all $n\in\mathbb{N}$. More precisely we have the following theorem.

\begin{theorem}[Staircase effect]\label{th:ROFstaircase}
Let $\lambda>4$, let $g_{n}$ be as in $(\ref{starcase
functions})$, and let $u_{n}$ be the unique minimizer of $\mathrm{ROF}%
_{\lambda,g_{n}}$ in $BV(]0,1[)$. Then for all $n$ sufficiently large there
exist $0<a_{n}<b_{n}<1$, with%
\[
a_{n}\rightarrow\frac{1}{\sqrt{\lambda}}\text{,\quad}b_{n}\rightarrow
1-\frac{1}{\sqrt{\lambda}}%
\]
as $n\rightarrow\infty$, such that $u_{n}=g_{n}$ on $\left[  a_{n}%
,b_{n}\right]  $ and $u_{n}$ is constant on each interval $\left[
0,a_{n}\right)  $ and $\left(  b_{n},1\right]  $.
\end{theorem}

\begin{proof}
Let $g_{n}^{-1}$ denote the generalized inverse function of $g_{n}$
defined by \eqref{gif}  with $g$ replaced by $g_{n}$. As both
$\{g_{n}\}$ and $\{g_{n}^{-1}\}$ converge uniformly to $g(x)=x$ and
since $\lambda>4$, one can check that for $n$ large enough there
exist $0<c_{1}^{\left(  n\right)  }<c_{2}^{\left( n\right)  }<1$
satisfying
\[
2\lambda\int_{0}^{g_{n}^{-1}(c_{1}^{\left(  n\right)  })}(c_{1}^{\left(
n\right)  }-g_{n})\,dx=1\qquad\text{and}\qquad2\lambda\int_{g_{n}^{-1}%
(c_{2}^{\left(  n\right)  })}^{1}(g_{n}-c_{2}^{\left(  n\right)  })\,dx=1\,
\]
with $c_{1}^{\left(  n\right)  }\rightarrow c_{1}$ and $c_{2}^{\left(
n\right)  }\rightarrow c_{2}$ as $n\rightarrow\infty$, where $c_{1}$ and
$c_{2}$ are defined by
\begin{equation}
2\lambda\int_{0}^{c_{1}}(c_{1}-x)\,dx=1\qquad\text{and}\qquad2\lambda
\int_{c_{2}}^{1}(x-c_{2})\,dx=1\,. \label{c1c2bis}%
\end{equation}
By Proposition~\ref{rof} the unique minimizer $u_{n}$ of $\mathrm{ROF}%
_{\lambda,g_{n}}$ in $BV(]0,1[)$ takes the form%
\[
u_{n}(x)=%
\begin{cases}
c_{1}^{\left(  n\right)  } & \text{if }0\leq x\leq g_{n}^{-1}(c_{1}^{\left(
n\right)  })\,,\\
g_{n}(x) & \text{if }g_{n}^{-1}(c_{1}^{\left(  n\right)  })<x\leq g_{n}%
^{-1}(c_{2}^{\left(  n\right)  })\,,\\
c_{2}^{\left(  n\right)  } & \text{if }g_{n}^{-1}(c_{2}^{\left(  n\right)
})<x\leq1\,.
\end{cases}
\]
The conclusion follows by observing that $a_{n}:=g_{n}^{-1}(c_{1}^{\left(
n\right)  })\rightarrow c_{1}$, $b_{n}:=g_{n}^{-1}(c_{2}^{\left(  n\right)
})\rightarrow c_{2}$ and that $c_{1}=\frac{1}{\sqrt{\lambda}}$ and
$c_{2}=1-\frac{1}{\sqrt{\lambda}}$, thanks to \eqref{c1c2bis}  .
\end{proof}

\subsection{Absence of the staircase effect: The case $p=1$}
Next we show that the presence of the higher order term in the functional
$\overline{\mathcal{F}}_1$ prevents the occurrence of the staircase effect. We
begin with the case $p=1$. We consider the minimization problem%

\begin{equation}
\min\Big\{  \overline{\mathcal{F}}_1\left(  u\right)  +\lambda\int_{a}%
^{b}(u-g)^{2}\,dx:\,u\in X^1_\psi({]a,b[})\Big\}
\,,\label{minimization p=1}%
\end{equation}
where $\overline{\mathcal{F}}_1$ is the relaxed functional given in
(\ref{relaxed}). To prove the absence of the staircase effect we need the
following auxiliary result that is of independent interest.

\begin{proposition}\label{Lipregp=1}
Assume that $\psi\colon{\mathbb{R}}\rightarrow{]0,+\infty\lbrack}$ is a
bounded Borel function satisfying $\left(  \ref{intfin}\right)  $ and  $\left(
\ref{inf>0}\right)$. Let $g:\left[  a,b\right]  \rightarrow{\mathbb{R}}$ be
Lipschitz continuous and let $u\in X^1_\psi({]a,b[})$ be a solution of the
minimization problem $(\ref{minimization p=1})$. Then $u$ is Lipschitz continuous and $u'\in BV({]a,b[})$.
\end{proposition}

\begin{proof}
The plan of the proof is the following. We will show that the
discontinuity set $S_u$ is empty and that  the  left and right
limits $(u')^a_-$ and $(u')^a_+$, defined in \eqref{uprime+}, are
finite everywhere on ${]a,b]}$ and on ${[a,b[}$, respectively. Note
that this will imply that the sets $Z^{\pm}[(u^{\prime})^{a}]$ (see
\eqref{Zu} and \eqref{Zu-}) are empty and, in turn, that $u\in
W^{1,1}({]a,b[})$ by the properties of the space
$X^1_\psi({]a,b[})$. Moreover, recalling that  the functions
$(u^{\prime})_{\vee}^{a}$  and $ (u^{\prime})_{\wedge}^{a}$ defined
in Remark~\ref{remark limits u'} are upper and lower semicontinuous
on $[a,b]$, it will also follow  that both $(u')^a_-$ and $(u')^a_+$
are bounded, yielding the Lipschitz continuity of $u$. In turn, the
fact that $u'\in BV({]a,b[})$ is a consequence of the local
Lipschitz continuity of $\Psi_1^{-1}$.

\noindent{\bf Step 1:} We start by showing that $S_u$ is empty.  We argue  by contradiction, assuming
that $S_u$ contains a point $x_0$. Without loss of generality we may suppose that $\nu_u(x_0)=1$; i.e.,
$u_+(x_0)> u_-(x_0)$.  We also assume that $\frac12 (u_+(x_0)+ u_-(x_0))\geq g(x_0)$.
In the following it is convenient to think of $u$ as coinciding everywhere with its lower semicontinuous representative
$u_{\wedge}:=\min\{u_-, u_+\}$.

Find $\varepsilon>0$ so small that
\begin{equation}\label{salti}
\sum_{\genfrac{}{}{0pt}2{\scriptstyle x\in S_u}{\scriptstyle x\in]x_0, x_0+\varepsilon[ }}
|[u](x)|<\frac{[u](x_0)}{4}
\end{equation}
and let $C>0$ satisfy
\begin{equation}\label{Csatisfy}
C>2\|g'\|_{\infty}\quad\text{and}\quad \frac12 (u_+(x_0)+ u_-(x_0))+C\varepsilon>u_-(x_0+\varepsilon)\,.
\end{equation}
For  $t\in [0,1]$ consider the affine function
$$
h^t(x):=\tfrac{(1-t)}{2}(u_+(x_0)+ u_-(x_0))+t\big(\tfrac14 u_-(x_0)+\tfrac34 u_+(x_0)\big)+C(x-x_0)
$$
and note that by \eqref{Csatisfy} there exists $x^t\in ]x_0, x_0+\varepsilon[$ such that
\begin{equation}\label{xt}
 (x^t,h^t(x^t))\in \Gamma_u\quad\text{and}\quad g<h^t<u\text{ in } ]x_0, x^t[\,,
\end{equation}
where $\Gamma_u$ stands for the extended graph of $u$ defined by
$$
\Gamma_u:=\{(x,t)\in {]a,b[}{\times}\R :
\min\{u_{-}(x),u_{+}(x)\}\le t\le \max\{u_{-}(x),u_{+}(x)\}\}\,.
$$
Let $u^t$ be the function defined by
\begin{equation}\label{ut}
u^t(x):=\begin{cases}
  h^t(x) & \text{if $x\in {]x_0, x^t[}$,}\\
u(x) & \text{otherwise,}
\end{cases}
\end{equation}
and note that
\begin{equation}\label{eta1}
\lambda\Big(\int_a^b| u-g|^2\, dx-\int_a^b| u^t-g|^2\, dx\Big)\geq
\lambda\Big(\int_a^b| u-g|^2\, dx-\int_a^b| u^1-g|^2\, dx\Big)=:\eta>0
\end{equation}
for every $t\in [0,1]$.
Now it is convenient to approximate $u$ with functions having only finitely many jump points.
 Hence the following approximation procedure is needed only when $S_u$ is infinite.
 In this case write $S_u=\{x_0, x_1, \dots, x_j, \dots\}$, for each $k$ define $S_{u}%
^{k}:=\left\{  x_{j}:\,0\leq j\leq k\right\}  $, and for $x\in{]a,b[}$ set
\[
u_{k}\left(  x\right)  =u_{+}\left(  a\right)  +\int_{a}^{x}\left(  u^{\prime
}\right)  ^{a}\,dt+\left(  u^{\prime}\right)  ^{c}\left(  \left]  a,x\right[
\right)  +\sum_{x_{j}<x,\,x_{j}\in S_{u}^{k}}\left[  u\right]  \left(
x_{j}\right)\,.
\]
Note that, since $u_k\to u$ in $L^\infty({]a,b[})$, for $k$ large enough
it follows from \eqref{Csatisfy} and \eqref{xt} that for every $t\in[0,1]$
there exists $x_k^t\in{]x_0, x_0+\varepsilon[}$ such that
$$
(x_k^t, h^t(x_k^t))\in \Gamma_{u_k}\quad\text{and}\quad g<h^t<u_k\text{ in } ]x_0, x_k^t[\,,
$$
where $ \Gamma_{u_k}$ denotes the extended graph of $u_k$.
 For all such $k$ we consider the comparison function $u_k^t$ defined as in \eqref{ut}, with $u$ and
 $x^t$  replaced by $u_k$ and $x^t_k$, respectively.
 Using the uniform convergence of $\{u_k\}$ to $u$ and \eqref{xt}, we have  that $x^t\le \liminf_k x^t_k$, which yields $u^t\ge \limsup_k u^t_k$ $\Lone$-a.e.\ on ${]a,b[}$.
 Moreover $u_k\to u$ in $\overline{\mathcal{F}}_1$ energy. Hence, also by \eqref{eta1}, we may find $k$ so large that for $t\in[0,1]$
 \begin{eqnarray}
& \displaystyle \!\!\!\!\!\!\!\!\!\!\lambda\Big(\int_a^b| u_k{-}g|^2\, dx-\int_a^b| u_k^t{-}g|^2\, dx\Big)\geq
\lambda\Big(\int_a^b| u_k{-}g|^2\, dx-\int_a^b| u_k^1{-}g|^2\, dx\Big)\geq \frac{\eta}{2}\,,
 \label{eta2}\\
& \displaystyle
 \overline {\mathcal F}_1( u_k)+\lambda\int_a^b| u_k-g|^2\, dx\leq\overline {\mathcal F}_1( u)+\lambda\int_a^b| u-g|^2\,
dx+\frac{\eta}{4}\,.\label{eta3}
 \end{eqnarray}
 Let us fix $k$ satisfying \eqref{eta2} and \eqref{eta3}.
 We claim that there exists $\bar t\in [0,1]$ such that $x_k^{\bar t}$ is a continuity point for $u_k$. Indeed, if not,
  then for every $t\in [0,1]$ there exists a jump point $x_j$, with $1\leq j\leq k$, such that $x_k^t=x_j$ and  the point $(x_k^t, h^t(x_k^t))$
  belongs to the corresponding vertical segment  of the extended graph of $u_k$.
  Setting $I_j:=\{t\in[0,1]:\, x_k^t=x_j\}$ and $\sigma_j:=\{(x_j, h^t(x_j)):\, t\in I_j\}$, it is clear that
  $[0,1]=\cup_{j=1}^k I_j$ and
  $\mathcal{H}^1(\sigma_j)=  \mathcal{H}^1(\{(x_0, h^t(x_0)):\, t\in I_j\})$. Thus,
 $$
 \sum_{\genfrac{}{}{0pt}2{\scriptstyle x\in S_u}{\scriptstyle x\in{]x_0, x_0+\varepsilon[ }}} |[u](x)|\geq
 \sum_{j=1}^k\mathcal{H}^1(\sigma_j)=\mathcal{H}^1\big(\{(x_0, h^t(x_0)):\, t\in [0,1]\}\big)=\frac{[u](x_0)}{4}\,,
$$
in contradiction with \eqref{salti}.

Since from now on $\bar t$ and $k$ are fixed, to simplify the notation we set $\hat x:=x_k^{\bar t}$, $\hat u:= u_k^{\bar t}$, $\hat h:= h^{\bar t}$, and
$\hat v:= \Psi_1\circ (\hat u')^a$.  By construction (see~\eqref{ut}) we have
\begin{equation}\label{varless}
|\hat u'|({]a,b[})\leq |u_k'|({]a,b[})\,.
\end{equation}
Next we claim that
\begin{equation}\label{nclaim}
( u')^a_-(\hat x)\leq\hat h'(\hat x)=C\,.
\end{equation}

If $( u')^a_-(\hat x)\leq 0$ there is nothing to prove. If  $( u')^a_-(\hat x)>0$, then by left continuity $( u')^a_-(y)>0$  for $y$ sufficiently close to $\hat x$, which, in turn, implies  $( u')^c(]y,\hat x[)\geq 0$
by the
 properties of $X^1_\psi({]a,b[})$.
 Since $S_{u_k}$ is finite and $\hat x$ is a continuity point, for $y$ in a left neighborhood of $\hat x$   we can write
\begin{equation*}
\hat h(\hat x)= u_k(\hat x)=u_k(y)+ \int_y^{\hat x}( u')^a(s)\, ds +(u')^c(]y,\hat x[)>
\hat h(y)+ \int_y^{\hat x}( u')^a(s)\, ds\,,
\end{equation*}
where we have used the fact that $u_k(\hat x)=\hat h(\hat x)$ and $\hat h< u_k$ in
 a left neighborhood of $\hat x$.
Claim \eqref{nclaim} follows.

Now, recalling that
$\Phi(1,t_1,t_2)=2\Psi_1(+\infty)-\Psi_1(t_1)-\Psi_1(t_2)$ for every
$t_1$, $t_2\in \overline\R$ by \eqref{Phi} and using
Remark~\ref{rem2.5'}, we estimate\begin{align}
|v^{\prime}|([x_{0},\hat{x}] &  \setminus S_{u})+\sum_{x\in
S_{u}\cap\lbrack
x_{0},\hat{x}]}\Phi(\nu_{u},(u^{\prime})_{-}^{a},(u^{\prime})_{+}%
^{a})\nonumber\\
\geq &  |v^{\prime}|({]x_{0},\hat{x}]})+\Phi(1,(u^{\prime})_{-}^{a}%
(x_{0}),(u^{\prime})_{+}^{a}(x_{0}))\nonumber\\
\geq &  |\Psi_{1}((u^{\prime})_{+}^{a}(x_{0}))-\Psi_{1}((u^{\prime})_{-}%
^{a}(\hat{x}))|+|\Psi_{1}((u^{\prime})_{+}^{a}(\hat{x}))-\Psi_{1}((u^{\prime
})_{-}^{a}(\hat{x}))|\nonumber\\
&  +\Phi(1,(u^{\prime})_{-}^{a}(x_{0}),(u^{\prime})_{+}^{a}(x_{0}))\nonumber\\
= &  |\Psi_{1}((u^{\prime})_{+}^{a}(x_{0}))-\Psi_{1}((u^{\prime})_{-}^{a}%
(\hat{x}))|+|\Psi_{1}((u^{\prime})_{+}^{a}(\hat{x}))-\Psi_{1}((u^{\prime}%
)_{-}^{a}(\hat{x}))|\nonumber\\
&  +2\Psi_{1}(+\infty)-\Psi_{1}((u^{\prime})_{-}^{a}(x_{0}))-\Psi
_{1}((u^{\prime})_{+}^{a}(x_{0}))\label{finest}\\
\geq &
-\Psi_{1}((u^{\prime})_{-}^{a}(\hat{x}))+2\Psi_{1}(+\infty)-\Psi
_{1}((u^{\prime})_{-}^{a}(x_{0}))+|\Psi_{1}((u^{\prime})_{+}^{a}(\hat
{x}))-\Psi_{1}((u^{\prime})_{-}^{a}(\hat{x}))|\nonumber\\
= &  \Psi_{1}(C)-\Psi_{1}((u^{\prime})_{-}^{a}(\hat{x}))+2\Psi_{1}%
(+\infty)-\Psi_{1}((u^{\prime})_{-}^{a}(x_{0}))-\Psi_{1}(C)\nonumber\\
&  +|\Psi_{1}((u^{\prime})_{+}^{a}(\hat{x}))-\Psi_{1}((u^{\prime})_{-}%
^{a}(\hat{x}))|\nonumber\\
\geq &  |\Psi_{1}(C)-\Psi_{1}((u^{\prime})_{+}^{a}(\hat{x}))|+\Phi
(1,(\hat{u}^{\prime})_{-}^{a}(x_{0}),(\hat{u}^{\prime})_{+}^{a}(x_{0}%
))\nonumber\\
= &  |\hat{v}^{\prime}|([x_{0},\hat{x}]\setminus
S_{\hat{u}})+\sum_{x\in S_{\hat{u}}\cap\lbrack
x_{0},\hat{x}]}\Phi(\nu_{\hat{u}},(\hat{u}^{\prime
})_{-}^{a},(\hat{u}^{\prime})_{+}^{a})\,,\nonumber
\end{align}
where in the last inequality we have used \eqref{ut} and~\eqref{nclaim}.
 Collecting \eqref{eta2},  \eqref{varless}, and \eqref{finest} we
 deduce that
 $$
 \overline {\mathcal F}_1(\hat u)+\lambda\int_a^b|\hat u-g|^2+\frac{\eta}{2}\leq \overline {\mathcal F}_1( u_k)+\lambda\int_a^b| u_k-g|^2
 $$
 and, in turn, by \eqref{eta3}
 \begin{equation}\label{claim1000}
 \overline {\mathcal F}_1(\hat u)+\lambda\int_a^b|\hat u-g|^2\, dx<\overline {\mathcal F}_1( u)+\lambda\int_a^b| u-g|^2\,
dx\,,
\end{equation}
which contradicts the minimality of $u$.

If $\frac12 (u_+(x_0)+ u_-(x_0))< g(x_0)$ then we proceed in a
similar manner: The comparison function  $\hat u$ is now constructed
by replacing $u_k$ with an affine function (defined as before and
with $C$ and $t$  properly chosen) in a left  neighborhood of $x_0$.
The argument is completely analogous to the previous one and we omit
the details.

\noindent{\bf Step 2:} We finally show that $(u')^a_-$ and $(u')^a_+$ are finite everywhere in ${]a,b]}$ and
in ${[a,b[}$, respectively. We give the details only for $(u')^a_-$, since one can argue for $(u')^a_+$ in an entirely
similar way.

Recall that by the previous step $u$ is continuous. Once again we reason by contradiction by assuming that there exists $\bar x\in {]a,b]}$ such that $|(u')^a_-(\bar x)|=+\infty$. Without loss
of generality we may suppose that $(u')^a_-(\bar x)=+\infty$. Using Remark~\ref{remark limits u'} and the differentiability
 properties of $BV$ functions we may choose a point $x_1\in{]a,\bar x[}$ such that $u$ is differentiable at $x_1$ and
\begin{equation}\label{choosex1}
u(x_1)\neq g(x_1)\,,\,\, u'(x_1)=(u')^a_-(x_1)=(u')^a_+(x_1)\,,\,\,
 u'(x_1)>2\|g'\|_{\infty},\,|v'|(]a,x_1])>0 \,.
\end{equation}
The first condition is a consequence of the fact that $g$ is Lipschitz and $u$ cannot be Lipschitz in any left neighborhood of $\bar x$, since $|(u')^a_-(\bar x)|=+\infty$.
The last condition follows easily from the fact that $(u')^a$ cannot be constant
$\Lone$-a.e. on
$]a, \bar x[$.
Assume that $u(x_1)>g(x_1)$. Then, by \eqref{choosex1} and by the previous step, we can find $\varepsilon\in {[0,\frac12[}$, with
\begin{equation}\label{withe}
\Psi_1(u'(x_1))-\Psi_1((1-\varepsilon)u'(x_1))<|v'|({[x_1,b[})\,,
\end{equation}
such that the affine function $h(x):=u(x_1)+(1-\varepsilon)u'(x_1)(x-x_1)$ satisfies one of the following
conditions: Either there exists a  point $x_2\in{ ]x_1,b[}$ for $u$  such that
\begin{equation}\label{hdiff}
 \text{$h(x_2)=u(x_2)$ and } g<h<u
\text{ in }]x_1,x_2[\,,
\end{equation}
or
\begin{equation}\label{or}
 g<h<u\text{ in }]x_1,b[\,.
\end{equation}
In the latter case we set $x_2:=b$.
We now consider the comparison function
$$
\hat u(x):=\begin{cases}
 h(x) & \text{if $x\in {]x_1,x_2[}$,}\\
u(x) & \text{otherwise,}
\end{cases}
$$
and we denote $\hat v:= \Psi_1\circ (\hat u')^a$.
We claim that  \eqref{claim1000} holds,
contradicting the minimality of~$u$.
By \eqref{hdiff} and \eqref{or} in any case we    have
$$
 \lambda\int_a^b|\hat u-g|^2\, dx<\lambda\int_a^b| u-g|^2\, dx\,.
$$
Moreover, if $x_2<b$ we have by construction
$|\hat u'|([x_1,x_2])=u(x_2)-u(x_1)\leq | u'|([x_1,x_2])$,
while if $x_2=b$ we have $|\hat u'|({[x_1,b[})=u_-(b)-u(x_1)\leq | u'|({[x_1,b[})$, so that
in both cases $|\hat u'|({[a,b[})\leq | u'|({]a,b[})$.
Hence \eqref{claim1000} will follow if we show that
$|\hat v'|([x_1,x_2])\leq | v'|([x_1, x_2])$, where $[x_1, x_2]$ is replaced by
${[x_1,b[}$ if $x_2=b$.
To see this we first assume that \eqref{hdiff} holds. Arguing as for \eqref{nclaim}, we deduce
$(u')^a_-(x_2)\leq h'(x_2)=(1-\varepsilon)u'(x_1)$. Therefore by \eqref{choosex1} we have
\begin{align*}
|v^{\prime}|([x_{1},x_{2}])= &
|v^{\prime}|({]x_{1},x_{2}[})+|v^{\prime
}|(\{x_{2}\})\\
\geq &  \Psi_{1}(u^{\prime}(x_{1}))-\Psi_{1}((u^{\prime})_{-}^{a}%
(x_{2}))+|\Psi_{1}((u^{\prime})_{-}^{a}(x_{2}))-\Psi_{1}((u^{\prime})_{+}%
^{a}(x_{2}))|\\
= &  \Psi_{1}(u^{\prime}(x_{1}))-\Psi_{1}((1-\varepsilon)u^{\prime}%
(x_{1}))+\Psi_{1}((1-\varepsilon)u^{\prime}(x_{1}))-\Psi_{1}((u^{\prime}%
)_{-}^{a}(x_{2}))\\
&  +|\Psi_{1}((u^{\prime})_{-}^{a}(x_{2}))-\Psi_{1}((u^{\prime})_{+}^{a}%
(x_{2}))|\\
\geq &  \Psi_{1}(u^{\prime}(x_{1}))-\Psi_{1}((1-\varepsilon)u^{\prime}%
(x_{1}))+|\Psi_{1}((1-\varepsilon)u^{\prime}(x_{1}))-\Psi_{1}((u^{\prime}%
)_{+}^{a}(x_{2}))|\\
= &  |\hat{v}^{\prime}|([x_{1},x_{2}])\,.
\end{align*}
If \eqref{or} holds then, by \eqref{withe}, we obtain
$$
  | v'|({[x_1, b[})  > \Psi_1(u'(x_1))- \Psi_1((1-\varepsilon)u'(x_1))=|\hat v'|([x_1,b[)\,.
$$
If $u(x_1)< g(x_1)$ we modify the previous argument in the following way. We now choose  $\varepsilon\in [0,\frac12[$ satisfying \eqref{withe} with $| v'|({[x_1, b[})$ replaced by
$| v'|({]a,x_1]})$
and such that the affine function $h(x)$ defined before  satisfies one of the following conditions: Either there exists a point $x_2\in {]a, x_1[}$ such that
$h(x_2)=u(x_2)$ and  $u<h<g$ in $]x_2,x_1[$,
or $u<h<g$ in $]a,x_1[$.
In the latter case we set $x_2:=a$.
We now consider the comparison function
$$
\hat u(x):=\begin{cases}
 h(x) & \text{if $x\in {]x_2,x_1[}$,}\\
u(x) & \text{otherwise,}
\end{cases}
$$
and we proceed exactly as before to show \eqref{claim1000}.
\end{proof}

We now turn to the main theorem of this subsection.

\begin{theorem}\label{th:m1}
Assume that $\psi\colon{\mathbb{R}}\rightarrow{]0,+\infty\lbrack}$
is a bounded Borel function satisfying \eqref{intfin} and
\eqref{inf>0}, let $g:\left[  a,b\right]  \rightarrow{\mathbb{R}}$
be Lipschitz continuous, and let $\{h_{n}\}$ satisfy
\begin{equation}\label{hto0}
h_n\rightharpoonup 0\quad\text{weakly$^*$ in }L^{\infty}\left(  \left]  a,b\right[  \right)\,.
\end{equation}
Define $\mathcal{A}_n$ as  the class of all solutions to  $(\ref{minimization p=1})$, with $g$ replaced by $g_{n}:=g+h_{n}$.
 Then for $n$ large enough  every solution $u_n\in\mathcal{A}_n$  is Lipschitz continuous. Moreover,
\begin{equation}\label{limsup1000}
 \limsup_{n\to\infty}\sup_{w\in \mathcal{A}_n}\|w\|_{1,\infty}<+\infty
\end{equation}
 and for every sequence $\{u_n\}\subset \mathcal{A}_n$ there exists a subsequence (not relabeled)  and
 a solution $u$ to \eqref{minimization p=1} such that $u_n\to u$ in $W^{1,p}(]a,b[)$ for all
 $p\in [1,+\infty[$.
\end{theorem}

\begin{proof}
It will be enough to prove that for any (sub)sequence
$\{u_n\}\subset\mathcal{A}_n$ we may extract a further subsequence
(not relabeled) and find a  solution $u$ to \eqref{minimization p=1}
such that $u_n$ is Lipschitz continuous for $n$ large enough,
\begin{equation}\label{showe}
 \limsup_{n\to\infty}\|u_n\|_{1,\infty}<+\infty\,,
\end{equation}
and $u_n\to u$  in $W^{1,p}(]a,b[)$ for all
 $p\in [1,+\infty[$.
Since the sequence $h_{n}$ is bounded in
$L^{\infty}\left(  \left]  a,b\right[  \right)  $ for any $w\in X^1_{\psi
}({]a,b[})$ we have%
\[
\sup_{n}\left(  \overline{\mathcal{F}}_1\left(  u_{n}\right)  +\lambda\int
_{a}^{b}(u_{n}-g_{n})^{2}\,dx\right)  \leq\overline{\mathcal{F}}_1\left(
w\right)  +\lambda\int_{a}^{b}(w-g_{n})^{2}\,dx\leq C<\infty\,,
\]
for a suitable constant $C>0$ independent of $n$.
By Corollary \ref{corollary compactness} there exist a subsequence not relabeled
and a function $u\in X^1_\psi({]a,b[})$ such that%
\begin{equation}
u_{n}\rightharpoonup u\quad\text{weakly}^{\ast}\text{ in }BV({]a,b[}%
)\,,\label{sc un weak star}%
\end{equation}
and%
\begin{equation}
u_{n}^{\prime}\rightarrow(u^{\prime})^{a}\quad\text{pointwise }\mathcal{L}%
^{1}\text{-a.e.\ in }{]a,b[}\,.\label{sc u'n pointwise}%
\end{equation}
Moreover, since also the functions   $h_{n}^{2}$  are equibounded,  upon extracting a further subsequence
 we may find $f\in L^{\infty}({]a,b[})$ such that
\begin{equation}\label{h2tof}
h_{n}^{2}\rightharpoonup f\quad\text{ weakly$^{\ast}$ in }L^{\infty}({]a,b[})\,.
\end{equation}
It is convenient to ``localize'' the functional $\overline{\mathcal{F}}_1$: For every Borel set $B\subset{]a,b[}$ and for $w\in X^1_\psi({]a,b[})$
we set
\begin{equation}
\overline{\mathcal{F}}_1(w;B):=|w^{\prime}|(B)+|v^{\prime}|(B\setminus
S_{w})+\sum_{x\in S_{w}\cap B}\Phi(\nu_{w},(w^{\prime})_{-}^{a},(w^{\prime
})_{+}^{a})\,,\label{local relaxed energy}%
\end{equation}
where $v:=\Psi_1\circ(w^{\prime})^{a}$.
We divide the remaining part the proof into two steps.

\noindent\textbf{Step 1:}
 We claim that $u$ is a
solution of the minimization problem $(\ref{minimization p=1})$ and that for
every open interval $I={]c,d[}$, with $a\leq c<d\leq b$ and $c$, $d\in [a,b]\setminus S_{(u')^a}$,
\begin{equation}
\lim_{n\rightarrow\infty}  \overline{\mathcal{F}}_1\left(  u_{n};I\right)
 =\overline{\mathcal{F}}_1\left(u;I\right)\,.\label{sc energy n goto}%
\end{equation}
To see this, note that for each $n\in\mathbb{N}$%
\[
\lambda\int_{I}(u_{n}-g_{n})^{2}\,dx=\lambda\int_{I}(u_{n}-g)^{2}\,dx-2\lambda\int
_{I}\left(  u_{n}-g\right)  h_{n}\,dx+\lambda\int_{I}h_{n}^{2}\,dx\,.
\]
By \eqref{hto0}, \eqref{sc un weak star},  and \eqref{h2tof} it follows that
\begin{equation}\label{continuo}
\lim_{n\rightarrow\infty}\int_{I}(u_{n}-g_{n})^{2}\,dx=\int_{I}(u-g)^{2}%
\,dx+\int_{I}f\,dx\,.
\end{equation}
Recall also that   by  lower semicontinuity
\begin{equation}
\liminf_{n\rightarrow\infty}\overline{\mathcal{F}}_1\left(
u_{n};A\right)   \geq
\overline{\mathcal{F}}_1\left(  u;A\right)\,,\label{sc lower ineq}%
\end{equation}
for every open set $A\subset{]a,b[}$.

By the minimality of $u_{n}$\ for every $w\in X^1_{\psi
}({]a,b[})$ we have%
\begin{align*}
  \overline{\mathcal{F}}_1\left(  w\right)  +&\lambda\int_{a}^{b}(w-g)^{2}%
\,dx-2\lambda\int_{a}^{b}\left(  w-g\right)  h_{n}\,dx+\lambda\int_{a}%
^{b}h_{n}^{2}\,dx\\
&  =\overline{\mathcal{F}}_1\left(  w\right)  +\lambda\int_{a}^{b}(w-g_{n}%
)^{2}\,dx\geq\overline{\mathcal{F}}_1\left(  u_{n}\right)  +\lambda\int_{a}%
^{b}(u_{n}-g_{n})^{2}\,dx\,.
\end{align*}
Using (\ref{sc lower ineq}) (with $A={]a,b[}$) and once again \eqref{hto0} and \eqref{h2tof}, we
get%
\begin{align*}
\overline{\mathcal{F}}_1&\left(  w\right)     +\lambda\int_{a}^{b}%
(w-g)^{2}\,dx+\lambda\int_{a}^{b}f\,dx\geq
 \limsup_{n\rightarrow\infty}\Big(
\overline{\mathcal{F}}_1\left(  u_{n}\right)  +\lambda\int_{a}^{b}(u_{n}%
-g_{n})^{2}\,dx\Big) \\
& \geq\liminf_{n\rightarrow\infty}\Big(
\overline{\mathcal{F}}_1\left(
u_{n}\right)  +\lambda\int_{a}^{b}(u_{n}-g_{n})^{2}\,dx\Big)  \geq \overline{\mathcal{F}}_1\left(  u\right)  +\lambda\int_{a}^{b}(u-g)^{2}%
\,dx+\lambda\int_{a}^{b}f\,dx\,.
\end{align*}
Given the arbitrariness of \ $w\in X^1_\psi({]a,b[})$ this implies
that $u$ is a solution of the minimization problem
$(\ref{minimization p=1})$. Moreover, taking $w=u$ in the previous
inequalities and using \eqref{continuo} we deduce (\ref{sc energy n
goto}) for $I={]a,b[}$; i.e.,
\begin{equation}\label{I={]a,b[}}
\lim_{n\to\infty}\overline {\mathcal {F}}_1(u_n)=\overline {\mathcal {F}}_1(u)\,.
\end{equation}

It remains to prove (\ref{sc energy n goto}) for every  open interval of the form
$I={]c,d[}$, with  $c$,$d\in [a,b]\setminus S_{(u')^a}$. To this end  fix one such interval and assume by contradiction that
\begin{equation}\label{assurdop=1}
\limsup_{n\rightarrow\infty} \overline{\mathcal{F}}_1\left(  u_{n};I\right)
  >\overline{\mathcal{F}}_1\left(
u;I\right) \,.
\end{equation}
As $u$ is continuous by Proposition~\ref{Lipregp=1}, our assumption on $I$ implies that the end points $c$ and $d$  do not charge  $ \overline{\mathcal{F}}_1\left(u;\cdot\right)$,   so that  $\overline{\mathcal{F}}_1\left(u;I\right)=\overline{\mathcal{F}}_1\left(
u;\overline I\cap{]a,b[}\right)$. Therefore, combining \eqref{sc lower ineq},
\eqref{I={]a,b[}}, and \eqref{assurdop=1}
 we obtain
\begin{eqnarray*}
&\displaystyle\overline {\mathcal {F}}_1(u)=\overline{\mathcal{F}}_1\left(
u;\overline I\cap{]a,b[}\right)+\overline{\mathcal{F}}_1\left(u;{]a,b[}\setminus\overline I\right)
=\overline{\mathcal{F}}_1\left(
u;I\right)+\overline{\mathcal{F}}_1\left(u;{]a,b[}\setminus\overline I\right)\vphantom{\int}\\
&\displaystyle\vphantom{\int}< \limsup_{n\rightarrow\infty} \overline{\mathcal{F}}_1\left(  u_{n};I\right)+
 \liminf_{n\rightarrow\infty} \overline{\mathcal{F}}_1\left(u_n;{]a,b[}\setminus\overline I\right)\leq
 \lim_{n\to\infty}\overline {\mathcal {F}}_1(u_n)=\overline {\mathcal {F}}_1(u)\,,
\end{eqnarray*}
which is a contradiction.
This concludes the proof of \eqref{sc energy n goto}.

\noindent\textbf{Step 2:} We now show that $u_n$ is Lipschitz continuous for $n$ large enough
and that \eqref{showe} holds. Note that the convergence of $u_n$ to $u$ in $W^{1,p}(]a,b[)$ for all $p\in[1,+\infty[$ will then easily follow
  from  \eqref{showe} and \eqref{sc u'n pointwise}. Assume by contradiction that the conclusion is false. Then, arguing as at the beginning of the proof of Proposition~\ref{Lipregp=1},
we may find a subsequence (not relabeled) and points
 $x_n\in {]a,b[}$ such that one of the following two
cases holds:
\begin{itemize}
 \item[(i)] $x_n\notin S_{(u'_n)^a}$ and $|(u'_n)^a(x_n)|\to+\infty$\smallskip;
\item[(ii)] $x_n\in S_{u_n}$  for every $n\in\N$.
\end{itemize}
Assume that (i) holds and, without loss of generality, that $(u'_n)^a(x_n)\to+\infty$. Upon extracting a
further subsequence we may also assume that $x_n\to x_0\in[a,b]$. Recall that by Proposition~\ref{Lipregp=1}
 and by the previous step the function $u$ is Lipschitz continuous. Hence there are two cases: Either
\begin{equation}\label{x00}
\overline{\mathcal{F}}_1(u;\{x_0\}\cap{]a,b[})=0
\end{equation}
or
\begin{equation}\label{x0>0}
x_0\in S_{(u')^a}\,,\quad\ (u')^a_{\pm}(x_0)\in \R\,, \quad\
 \overline{\mathcal{F}}_1(u;\{x_0\})=|\Psi_1((u')^a_{+}(x_0))-\Psi_1((u')^a_{-}(x_0))|\,.
\end{equation}
Assume first that \eqref{x00} holds.
 Set $L:=\|u'\|_{\infty}$ and fix $\e$
so small that,
$$
 \overline{\mathcal{F}}_1\left(u;I_\e\right)<\int_{L+1}^{+\infty}\psi(t)\, dt\,,
$$
where $I_\e:={]x_0-\e, x_0+\e[}\cap{]a,b[}$.
 By (\ref{sc energy n goto}) we also have
\begin{equation}\label{esmall}
 \overline{\mathcal{F}}_1\left(u_n;I_\e\right)<\int_{L+1}^{+\infty}\psi(t)\, dt\,,
\end{equation}
for $n$ large enough. On the other hand by \eqref{sc u'n pointwise} there exists $y\in I_\e$
such that $(u'_n)^a(y)<L+1$ for $n$ large. Moreover, taking into account (i), we also have
$(u'_n)^a(x_n)>L+1$  for $n$ large enough.
Thus,
$$
 \overline{\mathcal{F}}_1\left(u_n;I_\e\right)\geq |v'_n|(I_\e)\geq
 |\Psi_1((u'_n)^a(x_n))-\Psi_1((u'_n)^a(y))|\geq
 \Psi_1((u'_n)^a(x_n))-\Psi_1(L+1)\,,
$$
for all $n$ sufficiently large.
 Passing to the limit as $n\to\infty$ we then obtain
 $$
\liminf_{n\to\infty}  \overline{\mathcal{F}}_1\left(u_n;I_\e\right)\geq
\Psi_1(+\infty)-\Psi_1(L+1)= \int_{L+1}^{+\infty}\psi(t)\, dt\,,
$$
which contradicts \eqref{esmall}.

In case  \eqref{x0>0} holds, then $x_0\in{]a,b[}$. Set
\begin{equation}\label{seteta}
\eta:=2\Psi_1(+\infty)\!-\!\Psi_1((u')^a_+(x_0))\!-\!\Psi_1((u')^a_-(x_0))\!-
\!|\Psi_1((u')^a_{+}(x_0))\!-\!\Psi_1((u')^a_{-}(x_0))|>0
\end{equation}
and choose $\e$ such  that both $x_0-\e$ and  $x_0+\e$ belong to
${]a,b[}\setminus S_{(u')^a}$ and
\begin{eqnarray}
&  \displaystyle\overline{\mathcal{F}}_1\left(  u;I_\e\right)<
 |\Psi_1((u')^a_{+}(x_0))-\Psi_1((u')^a_{-}(x_0))|+\frac\eta3\,,
\label{etae}\\
&\displaystyle |\Psi_1((u')^a_{\pm}(y))-\Psi_1((u')^a_\pm(x_0))|
<\frac{\eta}{4}\quad \text{ for } y\in I^\pm_\e\,,\label{etae-bis}
\end{eqnarray}
where $I_\e:={]x_0-\e, x_0+\e[}$, $I^+_\e:={]x_0, x_0+\e[}$, and $I^-_\e:={]x_0-\e, x_0[}$.
Note that by \eqref{sc energy n goto} and \eqref{etae} we have
\begin{equation}\label{small1000}
 \overline{\mathcal{F}}_1\left(  u_{n};I_\e\right)< |\Psi_1((u')^a_{+}(x_0))-\Psi_1((u')^a_{-}(x_0))|+\frac\eta3
\end{equation}
for $n$ large enough. Moreover, by \eqref{sc u'n pointwise} and \eqref{etae-bis} we may find $y^-$,
$y^+\in I_\e$, with $y^-<x_0<y^+$, such that
\begin{equation} \label{etae2}
y^\pm\notin S_{(u'_n)^a}\quad\text{and}\quad |\Psi_1((u'_n)^a(y^\pm))-\Psi_1((u')^a_\pm(x_0))|<\frac{\eta}{4}
\end{equation}
 for $n$ large enough.
As $y^-<x_n<y^+$ for  $n$ sufficiently large, we  have
\begin{align}
 \overline{\mathcal{F}}_1\left(u_n;I_\e\right)\geq& |v'_n|(I_\e)\geq
 |\Psi_1((u'_n)^a(x_n))-\Psi_1((u'_n)^a(y^-))|\nonumber\\
 &+|\Psi_1((u'_n)^a(x_n))-\Psi_1((u'_n)^a(y^+))|\label{allsuch}\\
\geq&
 |\Psi_1((u'_n)^a(x_n))-\Psi_1((u')^a_-(x_0))|
 +|\Psi_1((u'_n)^a(x_n))-\Psi_1((u')^a_+(x_0))|-\frac\eta2\,,\nonumber
\end{align}
where the last inequality follows from \eqref{etae2}. Letting $n\to\infty$ in \eqref{allsuch} and
recalling \eqref{seteta} we deduce
\begin{eqnarray*}
&\displaystyle\liminf_{n\to\infty}  \overline{\mathcal{F}}_1\left(u_n;I_\e\right)\geq
2\Psi_1(+\infty)-\Psi_1((u')^a_+(x_0))-\Psi_1((u')^a_-(x_0))-\frac\eta2=\\
&\displaystyle =|\Psi_1((u')^a_{+}(x_0))-\Psi_1((u')^a_{-}(x_0))|+\frac\eta2\,,
\end{eqnarray*}
which contradicts \eqref{small1000}.
This concludes the proof of \eqref{showe} if (i) holds. An entirely similar argument can be used to treat the other  case.
\end{proof}

\subsection{Absence of the staircase effect: The case $p>1$}
We now turn to the case $p>1$. We consider the minimization problem%

\begin{equation}
\min\Big\{  \overline{\mathcal{F}}_{p}\left(  u\right)  +\lambda\int_{a}%
^{b}|u-g|^{2}\,dx:\,u\in X_{\psi}^{p}({]a,b[})\Big\}
\,,\label{minimization p>1}%
\end{equation}
where $\overline{\mathcal{F}}_{p}$ is the relaxed functional given in
(\ref{relaxed representation p>1}). We start with two auxiliary results.

\begin{proposition}
\label{prop:c1conv} Let $p>1$ and assume that $\psi\colon{\mathbb{R}%
}\rightarrow{]0,+\infty\lbrack}$ is a bounded Borel function satisfying
$\left(  \ref{inf>0}\right)  $ and $\left(  \ref{intfinp>1}\right)  $. Let $g$ be Lipschitz continuous  and let $u_n$ be a sequence in $X_{\psi
}^{p}({]a,b[})$ such that
$ \sup_n \overline{\mathcal{F}}_{p}(u_n)<+\infty$
and $u_n\to g$ in $L^2({]a,b[})$. Then $g\in C^1([a,b])\cap X_{\psi}^{p}({]a,b[})$. Moreover,  $u_n\in C^1([a,b])$ for $n$ large enough and $u_n\to g$ in $C^1([a,b])$.
\end{proposition}

\begin{proof}
By the assumptions and by Corollary~\ref{corollary compactness p>1} we deduce that
$g\in X_{\psi}^{p}({]a,b[})$. The fact that $g\in C^1([a,b])$ now follows from Remark~\ref{rm:definition}-(i).
To prove the last part of the statement  we start by showing that $(u_n^{\prime}%
)^{a}\rightarrow g^{\prime}$ uniformly in ${]a,b[}$.
Again by Corollary
\ref{corollary compactness p>1} the whole sequence $u_{n}$ satisfies
\begin{equation}
  \Psi_{p}\circ(u_{n}^{\prime})^{a}\rightharpoonup\Psi_{p}\circ
g^{\prime}\quad\text{weakly in }W^{1,p}({]a,b[})\,, \label{ulgotog}
\end{equation}
which implies, in particular, that
\begin{equation}\label{nle}
(\Psi_{p}\circ(u_{n}^{\prime})^{a})([a,b])\subset [ \Psi_{p}(-2\|g'\|_{\infty}), \Psi_{p}(2\|g'\|_{\infty}) ]\quad\text{for $n$ large enough.}
\end{equation}
Since by \eqref{inf>0}  $\Psi_{p}^{-1}$ is Lipschitz continuous on $[ \Psi_{p}(-2\|g'\|_{\infty}), \Psi_{p}(2\|g'\|_{\infty}) ]$, it follows from \eqref{ulgotog} and \eqref{nle} that $(u_n^{\prime}%
)^{a}\rightarrow g^{\prime}$ uniformly in ${]a,b[}$.
 In turn, by
Definition \ref{def spacep>1} we have that $u_{n}^{\prime}=(u_{n
}^{\prime})^{a}$ in $\left]  a,b\right[  $. In particular $u_{n}\in
C^{1}\left(  \left[  a,b\right]  \right)  $ by Remark~\ref{rm:definition}-(i)
and $u_{n}\rightarrow g$ in
$C^{1}\left(  \left[  a,b\right]  \right)$.
\end{proof}

\begin{proposition}
\label{prop:c1regp>1} Let $p$ and $\psi$ be as in the previous proposition. Then for every $C>0$ there exists
 $\overline\lambda=\overline\lambda(C)$ with the following property:  For all $g\in
C^{1}\left(  \left[  a,b\right]  \right)\cap X_{\psi}^{p}({]a,b[})$, with $\|g\|_{C^1([a,b])}\leq C$ and $\overline{\mathcal{F}}_{p}(g)\leq C$, and for all $\lambda\geq \overline\lambda$ every solution $u$ to \eqref{minimization p>1}  belongs to $C^1([a,b])$.
\end{proposition}
\begin{proof}
 Assume by contradiction that for every $n\in\mathbb{N}$ there exist $g_n\in
C^{1}\left(  \left[  a,b\right]  \right)\cap X_{\psi}^{p}({]a,b[})$, with $\|g_n'\|_{\infty}\leq C$ and $\overline{\mathcal{F}}_{p}(g_n)\leq C$, and a solution $u_n$ to
$$
\min\Big\{  \overline{\mathcal{F}}_{p}\left(  u\right)  +n\int_{a}%
^{b}|u-g_n|^{2}\,dx:\,u\in X_{\psi}^{p}({]a,b[})\Big\}
$$
which does not belong to $C^1([a,b])$. Owing to Proposition~\ref{prop:c1conv} we may assume, without loss of generality,  that
 $g_n\to g$ in $C^1([a,b])$ for a suitable function $g\in C^1([a,b])\cap X_{\psi}^{p}({]a,b[})$. Moreover, by
minimality, we have
$$
\overline{\mathcal{F}}_{p}\left(  u_n\right)  +n\int_{a}^{b}|u_n-g_n|^{2}\,dx\leq \overline{\mathcal{F}}_{p}(g_n)\leq C\,.
$$
It follows in particular that $ \sup_n \overline{\mathcal{F}}_{p}(u_n)<+\infty$ and $u_n\to g$ in $L^2({]a,b[})$. By Proposition~\ref{prop:c1conv} we conclude that $u_n\in C^1([a,b])$ for $n$ large enough, which gives a contradiction.
\end{proof}

The next theorem shows that also in the case $p>1$ the staircase effect does not occur.

\begin{theorem}
\label{theorem non staircase p>1} Let $\psi$ and $p$ be as in
 Proposition~\ref{prop:c1conv}, let $g\in
C^{1}\left(  \left[  a,b\right]  \right)\cap X_{\psi}^{p}({]a,b[}) $, and let $h_{n}$ satisfy \eqref{hto0}. For  $\lambda>0$ and $n\in\mathbb{N}$ let $\mathcal{A}_{\lambda,n}\subset X_{\psi}^{p}({]a,b[})$ be the class
 of the solutions to the minimization problem
$(\ref{minimization p>1})$, with $g$ replaced by $g_{n}:=g+h_{n}$. Let
$\overline{\lambda}$ be as in Proposition~\ref{prop:c1regp>1}, with
$C:=\max\{\|g\|_{C^1([a,b])},
 \overline{\mathcal{F}}_{p}(g)\}$. Then for all $\lambda\geq \overline{\lambda}$  we have
 $\mathcal{A}_{\lambda,n}\subset C^{1}\left(  \left[  a,b\right]  \right)$
 for $n$ sufficiently large. Moreover,
\begin{equation}\label{moreprecisely}
\lim_{\lambda\to\infty}\limsup_{n\to\infty}\sup_{u\in \mathcal{A}_{\lambda,n}}\| u- g\|_{C^1([a,b])}=0\,.
\end{equation}
\end{theorem}

\begin{proof}
We start by showing the second part of the statement.
Assume by contradiction that \eqref{moreprecisely} does not hold. Then there exist $\delta>0$, a sequence of real numbers
$\lambda_k\to +\infty$ and, for every $k$, a sequence of integers $n_j^k\to\infty$ as $j\to\infty$, such that for every $k$, $j$
\begin{equation}\label{deltaabs}
 \|u_{\lambda_k, n_j^k}-g\|_{C^1([a,b])}\geq \delta
\end{equation}
for a suitable function $u_{\lambda_k, n_j^k}\in \mathcal{A}_{\lambda_k, n_j^k}$ (with the
understanding that  $\|u_{\lambda_k, n_j^k}-g\|_{C^1([a,b])}=+\infty$ if
$u_{\lambda_k, n_j^k}\not\in C^1([a,b])$).
Arguing exactly as in Step 1 of the proof of Theorem~\ref{th:m1} we can show that for every $k$  there exist a subsequence (still denoted by $n_j^k$) and  a solution $u_{k}$  to \eqref{minimization p>1} with
$\lambda$ replaced by $\lambda_k$,
such that
\begin{equation}\label{eq1}
u_{\lambda_k, n_j^k}\rightharpoonup u_k\quad\text{weakly}^{\ast}\text{ in }BV({]a,b[})
\qquad\text{and}\qquad \overline{\mathcal{F}}_p(  u_{\lambda_k,n_j^k})\to \overline{\mathcal{F}}_p(  u_{k})
\end{equation}
as $j\to\infty$.  Moreover, since $g\in C^{1}\left(  \left[
a,b\right]  \right)\cap X_{\psi}^{p}({]a,b[}) $, we have by
minimality that
\begin{equation}\label{eq2}
 \overline{\mathcal{F}}_{p}\left(  u_{k}\right)  +\lambda_k\int_{a}%
^{b}|u_k-g|^{2}\,dx\leq  \overline{\mathcal{F}}_{p}(g)\,,
\end{equation}
which shows, in particular, that $u_{k}\to g$ in $L^2({]a,b[})$.
Combining \eqref{eq1} and \eqref{eq2},  and using a diagonal
argument, we may find a subsequence $n^k_{j_k}$ such that
$$
\sup_k \overline{\mathcal{F}}_{p}(u_{\lambda_k,n^k_{j_k}})<+\infty\qquad\text{and}\qquad
u_{\lambda_k,n^k_{j_k}}\to g\text{ in }L^2({]a,b[})\,.
$$
Proposition~\ref{prop:c1conv} then implies that $u_{\lambda_k,n^k_{j_k}}\to g$ in $C^1([a,b])$, which contradicts
\eqref{deltaabs}.

Finally, the first part of the statement follows from a similar argument by contradiction as a consequence of Propositions~\ref{prop:c1conv} and \ref{prop:c1regp>1} and of the fact that if
  $u_n\in\mathcal{A}_{\lambda,n}$ then,  up to subsequences, $u_n$ converges to a solution of \eqref{minimization p>1}.
\end{proof}

\noindent\textbf{Acknowledgments.} {The authors thank the Center for Nonlinear Analysis (NSF Grants No. DMS-0405343 and DMS-0635983) for its support during the preparation of this paper. The research of G.\ Dal Maso and M.\ Morini was partially supported by the Projects
``Calculus of Variations" 2004 and
``Problemi di Calcolo delle Variazioni in Meccanica e in Scienza dei Materiali'' 2006-2008, supported by the Italian Ministry of
Education, University, and Research, and by the  Project
``Variational Problems with Multiple Scales" 2006, supported by the Italian Ministry of
University and Research.
The research of I. Fonseca was partially supported by the National Science Foundation under Grant No. DMS-040171 and that of G. Leoni under Grants No. DMS-0405423 and DMS-0708039.
}

\newpage
\end{document}